\documentclass[11pt,a4paper]{article}
\usepackage[T1]{fontenc}

\usepackage{amsmath,amssymb,amsfonts,amsthm}

\usepackage[a4paper,margin=3cm]{geometry}
\usepackage{graphicx}
\usepackage{subcaption}
\usepackage{diagbox}
\usepackage{xcolor}

\usepackage{stmaryrd}

\newcommand{\Z}{\mathbb{Z}}
	
\newcommand{\N}{\mathbb{N}}	

\usepackage{booktabs}
\usepackage{tabularx}

\usepackage{enumitem}

\newtheorem{theorem}{Theorem}[section]
\newtheorem{lemma}[theorem]{Lemma}
\newtheorem{proposition}[theorem]{Proposition}

\theoremstyle{definition}
\newtheorem{definition}[theorem]{Definition}
\theoremstyle{remark}
\newtheorem{remark}[theorem]{Remark}

\usepackage[absolute,overlay]{textpos}

\usepackage{caption}
\usepackage{natbib}

\title{Construction, Enumeration and Duality of \\ Perfectly Balanced Musical Sets of Order $2pq$}

\author{Paul Lascabettes, Sidney Ali, Alexandre Popoff, Moreno Andreatta \\
\small CNRS, IRMA, ITI CREAA, Université de Strasbourg, France \\ 
\small STMS, IRCAM, CNRS, Sorbonne Université, Ministère de la Culture, Paris, France \\
}

\date{}

\begin{document}

\maketitle

\begin{abstract}
    In this paper, we present a comprehensive study of perfectly balanced musical sets, also known in the literature as vanishing sums of roots of unity, in the particular case of order $2pq$. Starting from the well-established result concerning the decomposition of such sets as an integer linear combination of regular polygons, we first provide an algebraic formalization that allows to explicitly construct all these sets. We then propose six different decompositions of each set into regular polygons. Moreover, by investigating their combinatorial enumeration, we derive a formula for the exact number of these sets for a given order $2pq$, and describe their distribution and cardinality for the particular case of order $6p$. One of the main results concerns the duality that can be established between these sets. In particular, we prove that every non-periodic minimal perfectly balanced set of order $2pq$ can be combined with another one to produce a self-complementary $2p$-periodic set, and the reciprocal also holds. Finally, we demonstrate that the modulus of the discrete Fourier transform coefficients of dual sets are related. All these results shed new light on the distribution of perfectly balanced sets of order $2pq$ and suggest new mathematical questions in a more general~case. 
\end{abstract}


\section{Introduction}

There are many possibilities for investigating the relationships between mathematics and music. For example, musical structures such as chords (\citet{reiner1985enumeration,fripertinger1991enumeration}), twelve-tone rows (\citet{read1997combinatorial}), musical mosaics (\citet{fripertinger1999enumeration}) or musical canons (\citet{fripertinger2001enumeration}) can all be described using an algebraic framework which offers to mathematicians the possibility to explore their combinatorial aspects. Music-theoretical problems can therefore be the starting point for advanced mathematical investigations, eventually leading to original constructions and new results in both fields. In this paper, we focus on musical objects which are subsets of the cyclic group $\Z_n$ with the property of being \textit{perfectly balanced}, defined as follows:

\begin{definition}
$A\subset\Z_n$ is \textit{perfectly balanced} if $\sum_{x \in A} e^{2i\pi x/n}=0$.
\end{definition}

This means that, by representing them as polygons inscribed into a circle, their centroid is located at the center of the circle. From a musical point of view, subsets of $\Z_n$ can be interpreted as musical scales or rhythms, and perfect balance corresponds to a distribution of notes with no particular attractor points. It is one of the geometrical properties dealing with the circular representation, such as maximally evenness (\citet{toussaint2005euclidean}), oddity (\citet{chemillier2003computation}) or several other properties that characterize musical sets (\citet{toussaint2019geometry,lascabettes2024good}).

Perfect balance has been originally introduced as a new way of constructing musical scales and rhythms (\citet{milne2015perfect}). Following this work, perfectly balanced sets have been described as an integer linear combination of regular polygons (\citet{milne2017exploring,amiot2017decompositions}), and further investigated within the framework of word theory (\citet{carey2017perfect}). More recently, the symmetric distribution of perfectly balanced sets of $\Z_n$ has been observed and characterized using the discrete Fourier transform (\cite{lascabettes2026perfectly}).
By analyzing the origin of perfectly balanced sets, it appears that this property was already known and studied by mathematicians as a particular case of \textit{vanishing sums of roots of unity} (see for example \citet{mann1965linear,conway1976trigonometric,poonen1998number,lam2000vanishing}). Although these works were developed in the framework of group rings, our paper focuses only on the particular case of elements of $\Z[\Z_n]$ with coefficients equal to 0 or 1 that we denote as \textit{sets}. 

We begin by providing some preliminary notations and definitions commonly used in the literature. We note $C_n$ the cyclic group of order $n$, $z$ a generator of $C_n$ and $\xi$ a primitive root of unity of order $n$. We consider the standard morphism $\varphi : \Z[C_n] \rightarrow \Z[\xi]$ defined by $\varphi(z)=\xi$ which enables to express perfectly balanced sets as the elements of $\ker(\varphi)$. We recall that for two elements $A=\sum_{i=0}^{n-1}a_iz^i$ and $B=\sum_{i=0}^{n-1}b_iz^i$ of $\Z[C_n]$, we have $z^j A=\sum_{i=0}^{n-1}a_iz^{i+j}$ and $A+B = \sum_{i=0}^{n-1}(a_i+b_i)z^i$. Moreover, the associated group ring is equipped with the order relation: $\sum_{i=0}^{n-1}a_iz^i \leq \sum_{i=0}^{n-1}b_iz^i$ iff $\forall i \in \llbracket 0, n-1 \rrbracket, a_i \leq b_i$. By an abuse of notation, we freely identify a subset $A\subset\Z_n$ with the corresponding element $A=\sum_{g\in\Z_n}a_i g \in \Z[C_n]$, where $a_g=1$ if $g\in A$ and $a_g=0$ otherwise.

Throughout this article, we focus on the case $n=2pq$, where $p$ and $q$ are always considered as distinct odd primes. Finally, we naturally restrict our study to \textit{non-periodic} \textit{minimal} perfectly balanced sets \textit{up to rotation}. They are defined as follows: 

\begin{definition}
Let $A\in \Z[C_n]$. The element $A$ is \textit{periodic} if there exists $t \in C_n$ with $t\neq 1$ such that $tA=A$. In this case, $A$ is said to be \textit{$t$-periodic} if $t$ is the smallest value satisfying the previous equality. It is \textit{non-periodic} otherwise.
\end{definition}

\begin{definition}
Let $A\in \Z[C_n]$ be a perfectly balanced set. The set $A$ is \textit{minimal} if for any perfectly balanced set $B$ such that $B \leq A$ either $B = 0$ or $B = A$.
\end{definition}

\begin{definition}
Let $A,B\in \Z[C_n]$. The sets $A$ and $B$ are \textit{equivalent up to rotation} if there exists $t\in C_n$ such that $tA=B$.
\end{definition}

To simplify the notations, we define $\Gamma_n$ as the class of all non-periodic minimal perfectly balanced sets of order~$n$. As already known in the literature, if $n$ admits only two distinct prime factors $p$ and $q$, i.e. $n=p^aq^b$, then the only minimal perfectly balanced sets of order~$n$ are periodic and they correspond to $C_p$ and $C_q$ up to rotation, which means that~$\Gamma_{p^aq^b}$ is empty (see Corollary~3.4, p.~98, \citet{lam2000vanishing}). Therefore, a necessary condition for the existence of non-periodic minimal perfectly balanced sets of $\Z_n$ is that $n$ admits at least three distinct prime factors. Moreover, in this case, the minimal perfectly balanced sets of order $p^aq^br^c$ are the same as those of order $pqr$ (see Theorem~3.1, p.~96, \citet{lam2000vanishing}). For this reason, by investigating the minimal perfectly balanced sets of order $2pq$, the case of order $2^a p^b q^c$ follows automatically. 

In this paper, we provide a detailed study of $\Gamma_{2pq}$ including an explicit construction, enumeration and classification of its elements. We reveal the symmetric character of their distribution with respect to the cardinality, as illustrated in Figure~\ref{fig:N42} for $\Gamma_{42}$ where $p=3$ and $q=7$. In this case, there are respectively 1, 3, 5, 5, 3, and 1 elements with cardinalities ranging from 8 to 13. In order to understand this symmetric distribution, we introduce a \textit{duality principle} that reveals the pairwise associations between all the elements of $\Gamma_{2pq}$. In particular, we prove that each dual association leads to a specific periodic and self-complementary element of $\Z[C_{2pq}]$. This duality principle adds an original property to the characterization of perfectly balanced sets.

\newcommand{\Nimg}[1]{\includegraphics[width=0.8\linewidth]{Figures/N_42/N42_#1.pdf}\par}
\begin{figure}[h]
    \centering
    \begin{subfigure}[c]{0.13\textwidth}
        \centering
        \Nimg{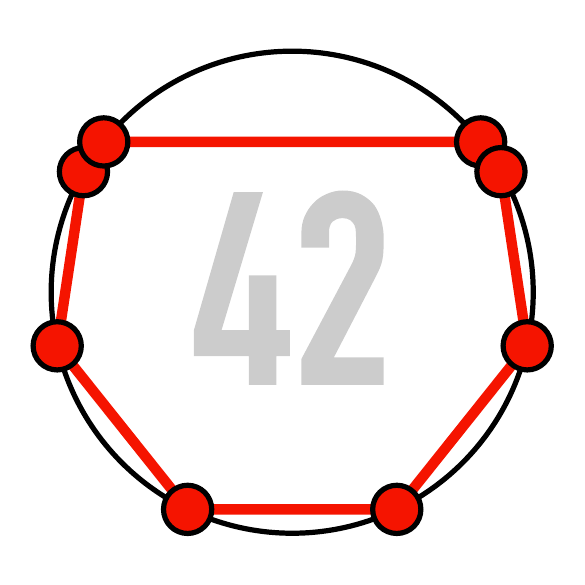}
        \caption{$k=8$}
    \end{subfigure}%
    \hspace{0.01\textwidth}
    \begin{subfigure}[c]{0.13\textwidth}
        \centering
        \Nimg{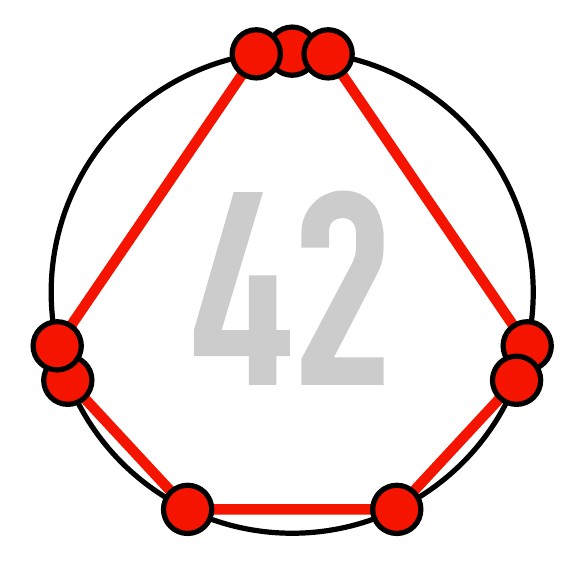}\Nimg{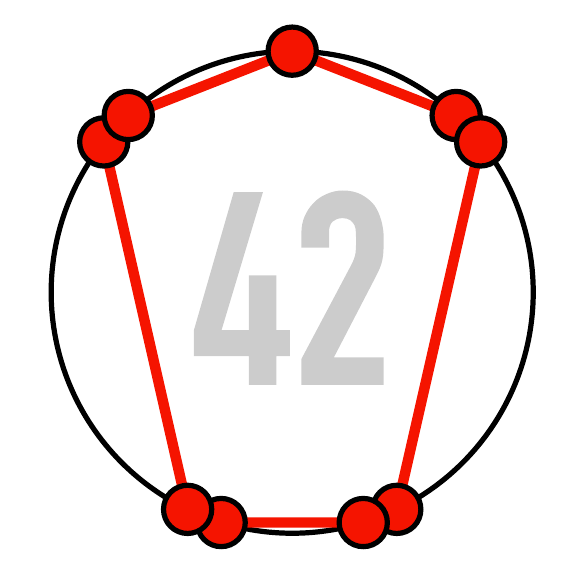}\Nimg{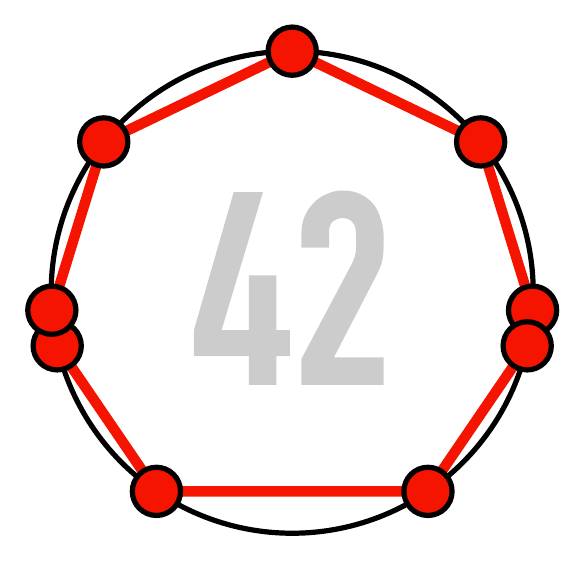}
        \caption{$k=9$}
    \end{subfigure}%
    \hspace{0.01\textwidth}
    \begin{subfigure}[c]{0.13\textwidth}
        \centering
        \Nimg{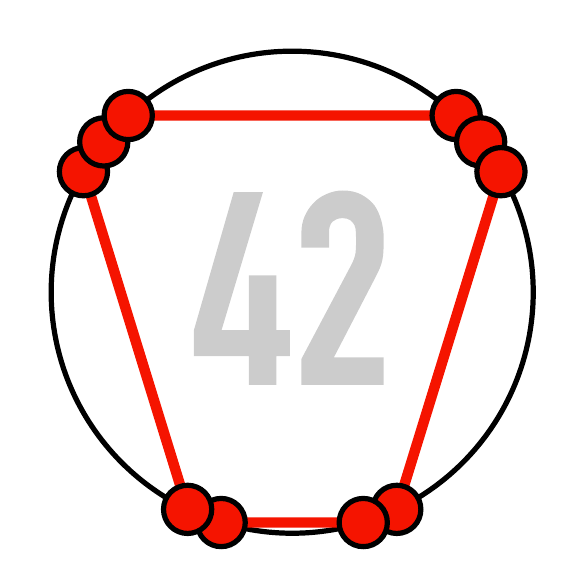}\Nimg{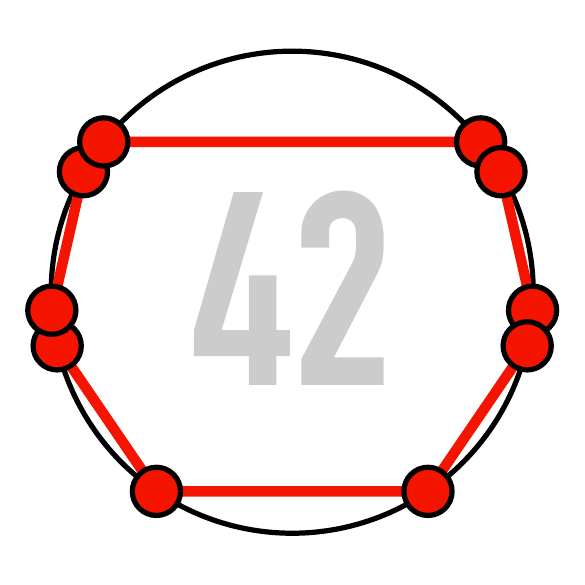}\Nimg{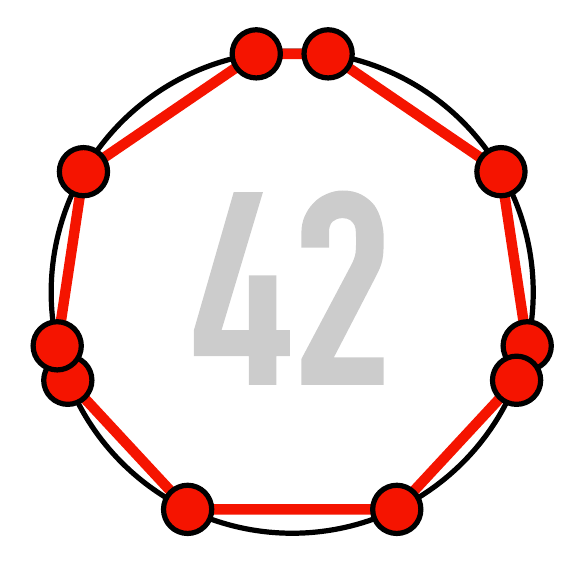}\Nimg{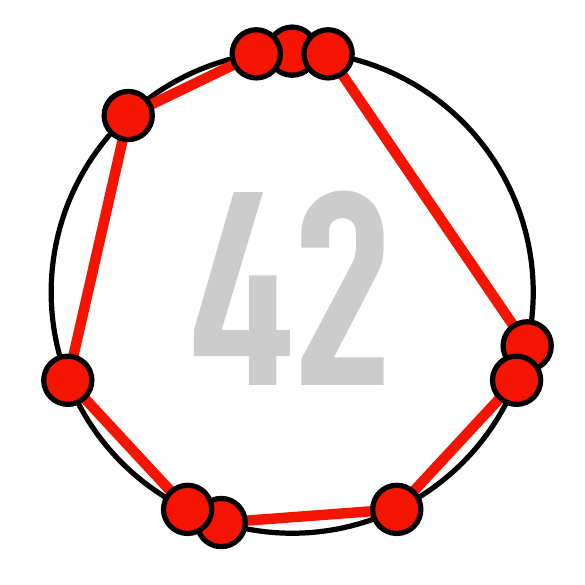}\Nimg{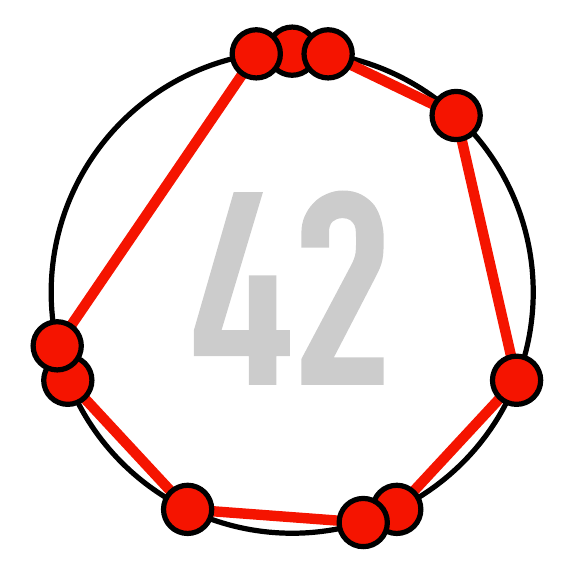}
        \caption{$k=10$}
        \label{fig:N42k10}
    \end{subfigure}%
    \hspace{0.01\textwidth}
    \begin{subfigure}[c]{0.13\textwidth}
        \centering
        \Nimg{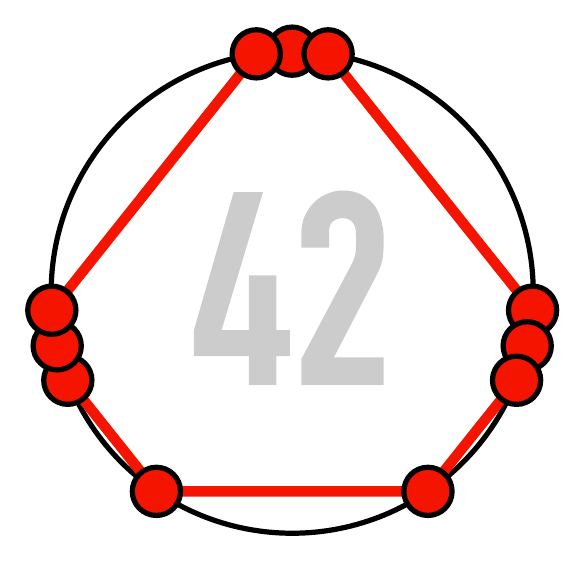}\Nimg{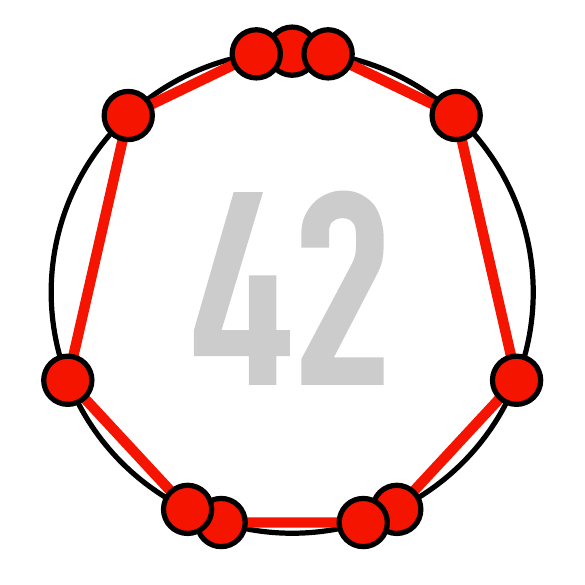}\Nimg{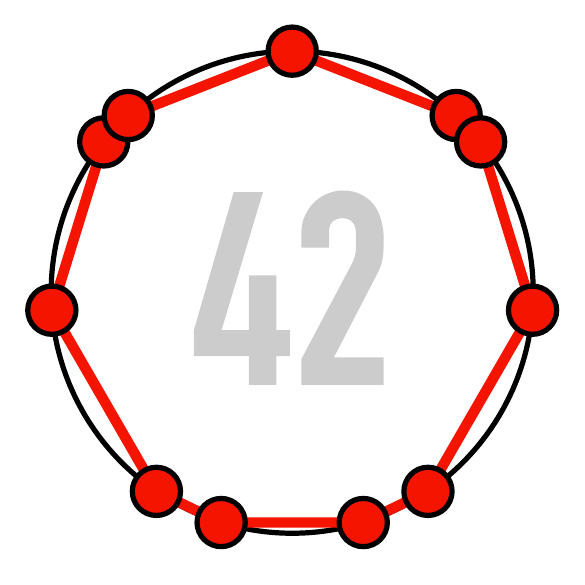}\Nimg{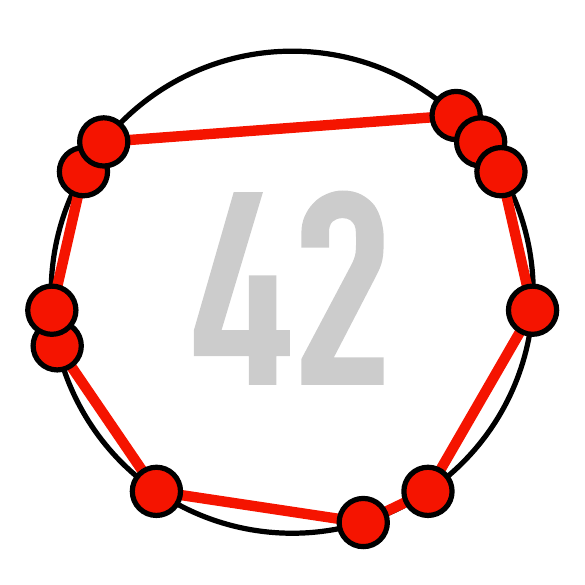}\Nimg{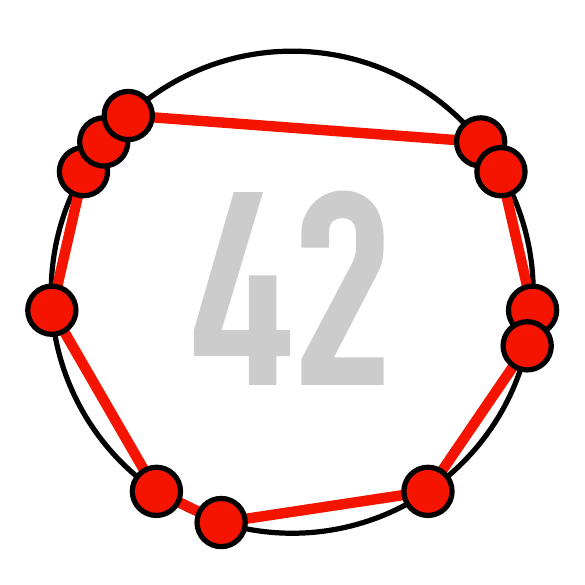}
        \caption{$k=11$}
        \label{fig:N42k11}
    \end{subfigure}%
    \hspace{0.01\textwidth}
    \begin{subfigure}[c]{0.13\textwidth}
        \centering
        \Nimg{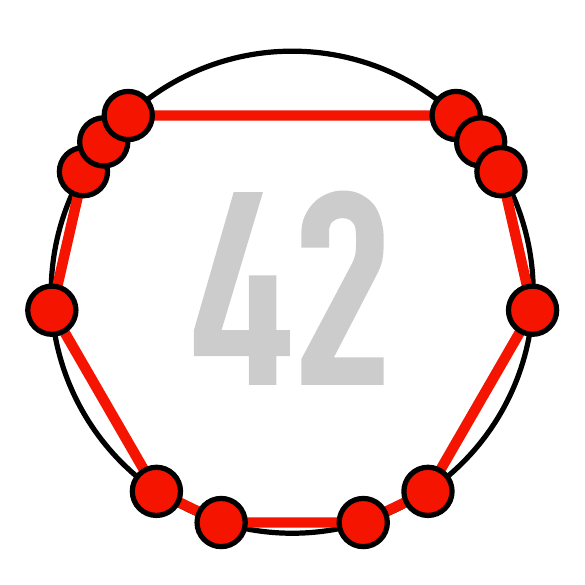}\Nimg{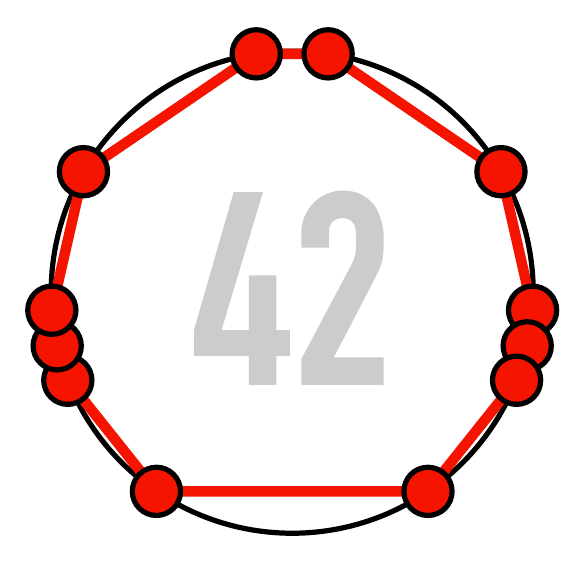}\Nimg{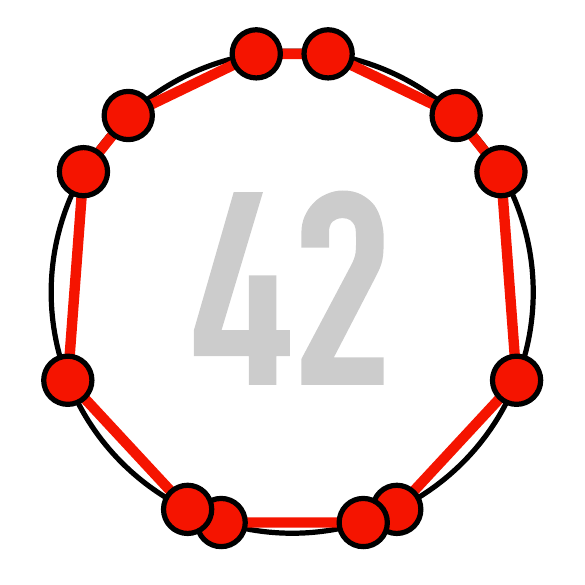}
        \caption{$k=12$}
    \end{subfigure}%
    \hspace{0.01\textwidth}
    \begin{subfigure}[c]{0.13\textwidth}
        \centering
        \Nimg{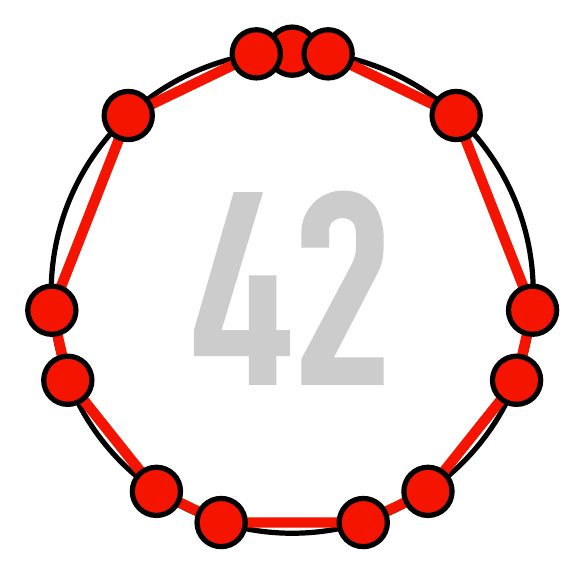}
        \caption{$k=13$}
    \end{subfigure}
    \caption{The 18 non-periodic minimal perfectly balanced sets of $\Z_{42}$ up to rotation and ordered by their cardinality~$k$.}
    \label{fig:N42}
\end{figure}

The rest of the paper is structured as follows. Section~\ref{sec:Algebraic} presents the algebraic construction of a perfectly balanced set of order $2pq$ (Section~\ref{subsec:Algebraic_Construction}) and shows six possible decompositions into regular polygons (Section~\ref{subsec:Algebraic_Decomposition}). Section~\ref{sec:Combinatorial} provides a combinatorial enumeration by establishing a formula for the number of non-periodic minimal perfectly balanced sets of order~$2pq$ up to rotation (Section~\ref{subsec:Combinatorial_Number}) and gives their distribution and cardinality for the particular case of order $6p$ (Section~\ref{subsec:Combinatorial_Distribution}). Section~\ref{sec:Pairwise} reveals the links between non-periodic minimal perfectly balanced sets of order $2pq$ by showing that a dual association produces a $2p$-periodic self-complementary set (Section~\ref{subsec:Pairwise_Main}). We prove that the reciprocal also holds, that is, every $2p$-periodic self-complementary set is the result of a dual association of two perfectly balanced sets, and we explain this association in the case of order $6p$ (Section~\ref{subsec:Reciprocal}). Section~\ref{sec:Fourier} introduces the relation between perfectly balanced sets and the discrete Fourier transform (Section~\ref{subsec:Fourier_Interpretation}) and proves that dual sets have almost the same modulus of the discrete Fourier transform coefficients (Section~\ref{subsec:Fourier_Relationship}). Finally, Section~\ref{sec:Conclusion} suggests some future research directions.


\section{Algebraic Construction and Decomposition of Perfectly Balanced Sets} \label{sec:Algebraic}

In this section, we present a way of describing non-periodic minimal perfectly balanced sets in the particular case of order $2pq$. The main result is a theorem that constructs such sets as a sum of two components which are given by disjoint unions of particular subgroups of prime order of $\Z_{2pq}$. This formula leads to six different decompositions of any non-periodic minimal perfectly balanced set of order $2pq$ into subgroups of prime order of $\Z_{2pq}$ that can be interpreted as regular polygons.


\subsection{Construction of Minimal Perfectly Balanced Sets} \label{subsec:Algebraic_Construction}

It is known in the literature that minimal perfectly balanced sets can be characterized as integer linear combinations of regular polygons. This result comes from the Rédei--De Bruijn--Schoenberg theorem (\citet{redei1950beitrag,debruijn1953factorization,schoenberg1964note}) that can be stated as follows, where $C_{p_i}$ can be interpreted as a regular $p_i$-gon and $\lambda_i C_{p_i}$ as a linear combination of its rotations.

\begin{theorem} \label{thm:redei-debuijn-schoenberg}
Let $A\in \Z[C_n]$ where $n=p_1^{\alpha_1}p_2^{\alpha_2} \cdots p_r^{\alpha_r}$ with $p_i$ distinct primes. Then:
$$ A \textit{ is perfectly balanced } \iff A = \sum_{i=1}^r \lambda_i C_{p_i}, \textit{ with } \lambda_i \in \Z[C_n].$$
\end{theorem}

H.~Mann successfully described linear relations between roots of unity with rational integer coefficients in the particular case of order $2pq$ (Theorem 4, p.~111, \citet{mann1965linear})). However, by focusing on sets, we can adapt this result and provide a general expression of non-periodic minimal perfectly balanced sets of order $2pq$ which allows to easily compute them all. To state this result, we specify that a partition $(P_1,P_2)$ of $C_p$ is said to be \textit{nontrivial} if both $P_1$ and $P_2$ are nonempty.

\begin{theorem} \label{thm:cas2pq}
Let $A\in \Z[C_n]$ with $n=2pq$, and let $\delta$ be the element of order $2$ of $C_n$. Then, $A$ is a non-periodic minimal perfectly balanced set if and only if there exist nontrivial partitions $(P_1,P_2)$ and $(Q_1,Q_2)$ of $C_p$ and $C_q$ such that:
$$A = \sum_{g \in P_1, h \in Q_1} gh + \sum_{g \in P_2, h \in Q_2} \delta gh.$$
\end{theorem}

\begin{remark} \label{rem:Notation_A_P1Q1}
To simplify the notation used in Theorem~\ref{thm:cas2pq}, we also write:
$$ A = P_1 Q_1 + \delta P_2 Q_2.$$
\end{remark}

\begin{proof}
Let $A \in \Gamma_{2pq}$.
Then $A$ gives a $\mathbb{Z}$-linear relation between $pq$-th roots of unity such that we can apply a result established by H. Mann (Theorem 4, p.~111, \citet{mann1965linear}).
Therefore, by writing $\rho ^{p} = \eta ^{q} = \delta^2 = 1$ with $\rho , \eta, \delta$ not equal to 1, $A$ is equal to one of the following expressions:
\begin{equation}\label{expression_1}
\sum_{i=1}^{t} \rho ^{\alpha_i} + \eta  \Sigma_1 + \cdots + \eta ^{q-1} \Sigma_{q-1},
\end{equation}
\begin{equation}\label{expression_2}
\delta\left(  \sum_{i=1}^{t} \rho ^{\alpha_i} + \eta  \Sigma_1 + \cdots + \eta ^{q-1} \Sigma_{q-1} \right),
\end{equation}
where $0 \leqslant \alpha_i \leqslant p-1$ for $ 1 \leqslant i \leqslant t$, $ 1 \leqslant t \leqslant p-1 $
and for $1 \leqslant j \leqslant q-1$, $\Sigma_j$ is either equal to:
\[ \sum_{i=1}^{t}  \rho ^{\alpha_i} \qquad \text{or} \qquad 
 \sum_{\substack{ \beta \neq  \alpha_i \\ 0 \leqslant  \beta \leqslant  p-1 } } \delta \rho ^{\beta},\] 
with this last expression appearing at least once. To get the announced expression, let $P_1 = \left\lbrace \rho ^{\alpha_i} \mid i \in  \llbracket 1, t \rrbracket \right\rbrace $ 
and set $P_2 = C_p \setminus P_1$.
Let also $Q_1$ be the set of $\eta^{j} \in C_q  $ such that $\eta^j \sum_{g \in P_1}^{} g \leqslant A$
if $A$ is equal to (\ref{expression_1}), to the set of $\eta ^{j} \in  C_{q}$ such that
$\delta \eta^{j} \sum_{g \in  P_1}^{} g \leqslant  A $ otherwise.
Set like before $Q_2 = C_q \setminus Q_1$.
Thus by construction, the pairs $(P_1, P_2)$ and $(Q_1, Q_2)$ are nontrivial partitions of $C_p$ and $C_q$ respectively.

If $A$ is equal to (\ref{expression_1}), we have as announced:
\[A = \sum_{g \in Q_1}^{}g \sum_{i=1}^{t}\rho ^{\alpha_i} + \sum_{g \in  Q_2}^{} g \sum_{\substack{ \beta \neq  \alpha_i \\ 0 \leqslant  \beta \leqslant  p-1 } }^{} \delta \rho ^{\beta}
= \sum_{g \in  Q_1}^{}\sum_{h \in  P_1} gh + \sum_{g \in  Q_2}^{} \sum_{h \in P_2 } \delta gh.\] 
If $A$ is equal to (\ref{expression_2}), a similar computation and a suitable relabelling of $P_1, P_2, Q_1$ and $Q_2$ yield the result.

Conversely, let $(P_1,P_2)$ and $(Q_1,Q_2)$ be nontrivial partitions of $C_p$ and $C_q$ such that
$A = P_1 Q_1 + \delta P_2 Q_2$, and let us show that $A$ is a non-periodic minimal perfectly balanced set.

Firstly, we show that $A$ is a perfectly balanced set. Since $C_q $ is perfectly balanced we have $\varphi(C_q) = 0$ which leads to $\varphi(Q_1) = - \varphi (Q_2) $. Therefore:
$$\varphi(A) = \varphi(P_1 Q_1) - \varphi(P_2 Q_2) = \varphi(P_1 Q_1) + \varphi(P_2 Q_1) = \varphi(C_p Q_1) = 0.$$

Secondly, we show that $A$ is non-periodic. To do that, it is enough to show that $A$ does not contain $C_2$, $C_p$ or $C_q$ up to rotation. Let $x \in A$, then there are two possibilities: either $x =g_1h_1$ with $g_1 \in P_1$ and $h_1 \in Q_1$ and $\delta x = \delta g_1h_1 \notin A$, or $x = \delta g_2h_2$ with $g_2 \in P_2$ and $h_2 \in Q_2$ and $\delta x = g_2h_2 \notin A$. In both cases $x + \delta x$ does not belong to $A$ which shows that $A$ does not contain $C_2$ up to rotation. Moreover, $A$ does not contain $C_p$ (respectively $C_q$) up to rotation because $P_1, P_2, Q_1$ and $Q_2$ are not empty.

Thirdly, we show that $A$ is minimal. Let $A'\leq A$ be a nonempty minimal perfectly balanced set. Then $A'$ is non-periodic since $A$ does not contain a nonempty periodic subset. Therefore, $A'$ is a minimal non-periodic perfectly balanced set. By applying the first implication of this theorem, we have $A'=P_1'Q_1'+\delta P_2'Q_2'$ with $(P_1',P_2')$ and $(Q_1',Q_2')$ nontrivial partitions of $C_p$ and $C_q$. Since $\Z_{2pq}\simeq \Z_2\times \Z_p\times \Z_q$, the inequality $A'\leq A$ implies $P_1'Q_1'\leq P_1Q_1$ and $P_2'Q_2'\leq P_2Q_2$. Therefore, $P_i'\leq P_i$ and $Q_i'\leq Q_i$ for $i=1,2$. Since $(P_1',P_2')$ and $(P_1,P_2)$ are partitions of $C_p$, we have $P_1'=P_1$ and $P_2'=P_2$. The same argument gives $Q_1'=Q_1$ and $Q_2'=Q_2$, which leads to the equality $A=A'$.
\end{proof}

Theorem~\ref{thm:cas2pq} provides a very efficient way to construct $\Gamma_{2pq}$ by considering the different nontrivial partitions $(P_1,P_2)$ and $(Q_1,Q_2)$ of $C_p$ and $C_q$ respectively. For example, in the case of $n=2\times 3 \times 5 = 30,$ there are three periodic minimal perfectly balanced sets, which are $C_2$, $C_3$ and $C_5$, and six non-periodic ones that are represented in Figure~\ref{fig:N30} using the circular representation. 

\newcommand{\NimgThirty}[1]{%
    \includegraphics[width=0.85\linewidth]{Figures/N_30/N30_#1.pdf}\par
}
\begin{figure}[h]
    \centering
    \begin{subfigure}[c]{0.22\textwidth}
        \centering
        \NimgThirty{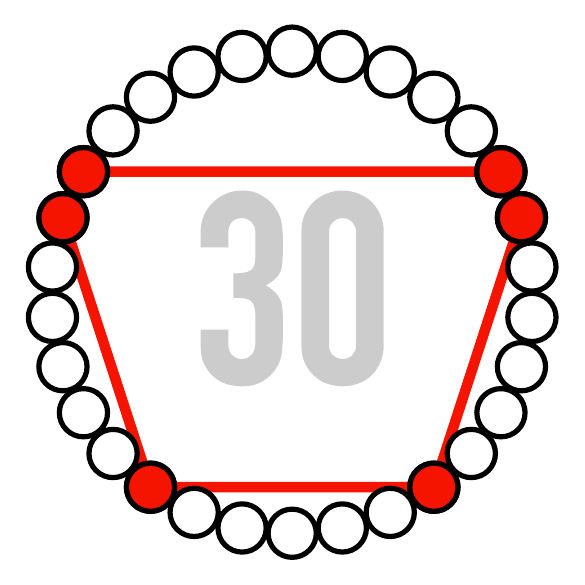}
        \caption{$k=6$}
        \label{subfig:N30_1}
    \end{subfigure}%
    \hspace{0.01\textwidth}
    \begin{subfigure}[c]{0.22\textwidth}
        \centering
        \NimgThirty{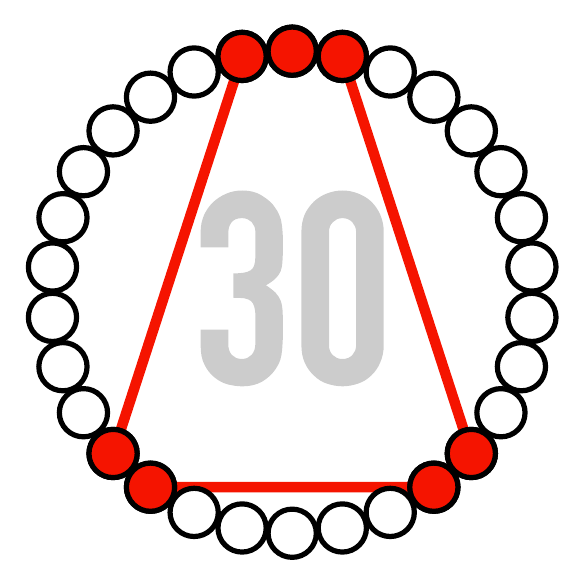} \vspace{0.4cm} \NimgThirty{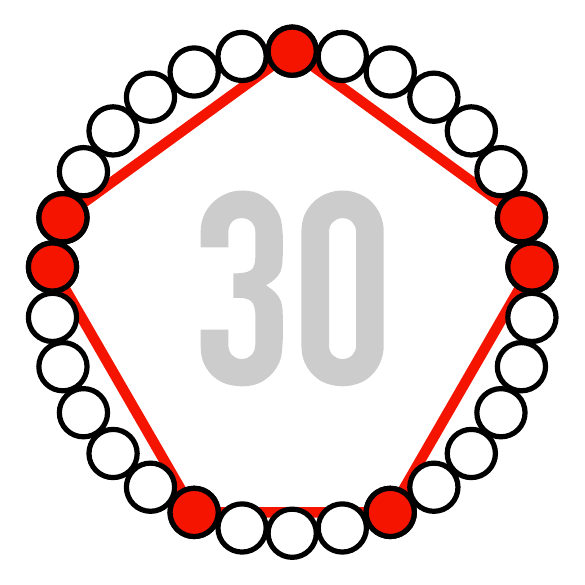}
        \caption{$k=7$}
        \label{subfig:N30_2}
    \end{subfigure}%
    \hspace{0.01\textwidth}
    \begin{subfigure}[c]{0.22\textwidth}
        \centering
        \NimgThirty{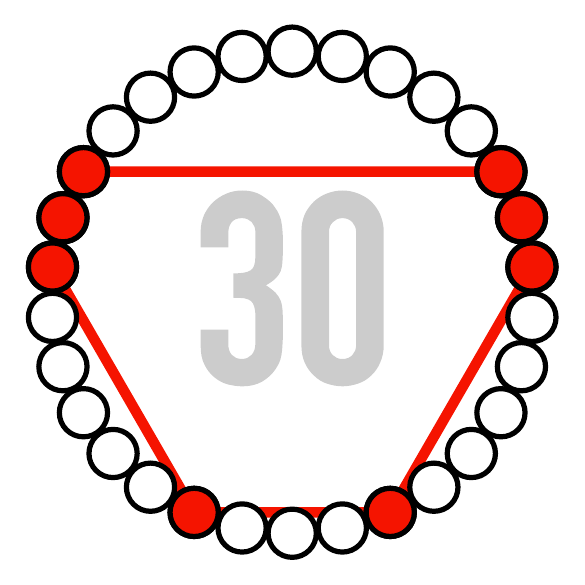} \vspace{0.4cm} \NimgThirty{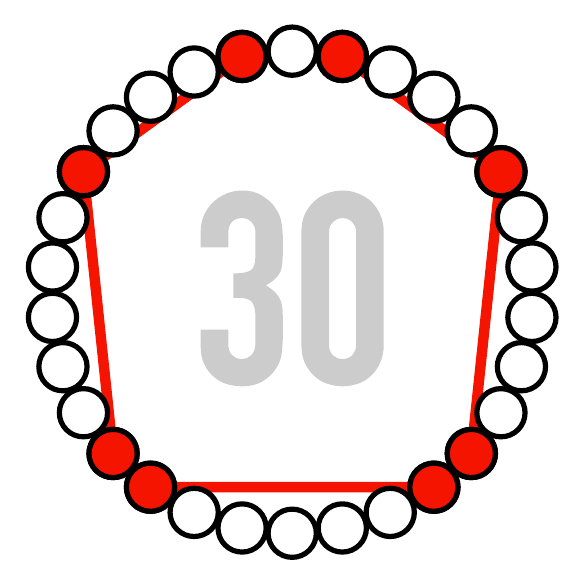}
        \caption{$k=8$}
        \label{subfig:N30_3}
    \end{subfigure}%
    \hspace{0.01\textwidth}
    \begin{subfigure}[c]{0.22\textwidth}
        \centering
        \NimgThirty{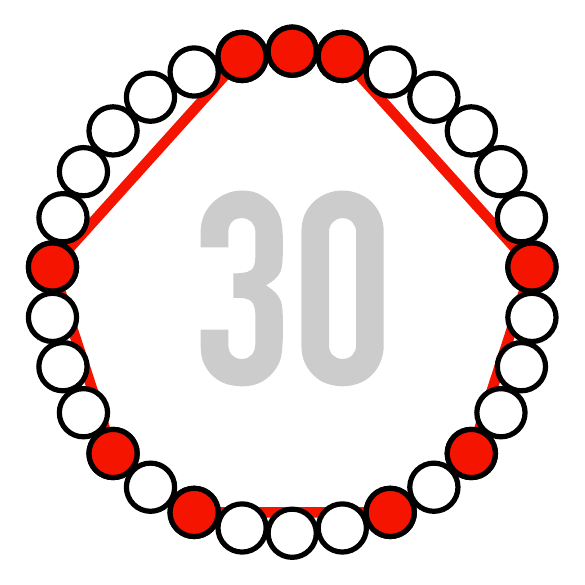}
        \caption{$k=9$}
        \label{subfig:N30_4}
    \end{subfigure}
    \caption{The six non-periodic minimal perfectly balanced sets of $\Z_{30}$, up to rotation, with cardinality $k$, illustrated using the circular representation.}
    \label{fig:N30}
\end{figure}

Their cardinalities are respectively equal to $k=6$, 7, 8 and 9. There have been many attempts to enumerate the way in which $k$ roots of unity can sum to zero. The first classification was provided for $k\leq 7$ (Theorem~6, p.~114, \citet{mann1965linear}), extended to $k\leq 9$ (Theorem~6, p.~237, \citet{conway1976trigonometric}), followed by $k \leq 12$ (Theorem~3.1, p.~141, \citet{poonen1998number}) and finally to $k \leq 16$ (Theorem~3.3, p.~8, \citet{christie2020classifying}). It is worth noting that all these minimal vanishing sums roots of unity with $k \leq 16$ have all their coefficients equal to 1 which mean that they can be interpreted as minimal perfectly balanced sets. It has been conjectured that the value of $k$ must be at least equal to 21 in order to have a minimal vanishing sums roots of unity with at least a coefficient higher than 1 (Conjecture 4.5, p.~16, \citet{christie2020classifying}). In that case, the value of the coefficients are 1 and 2, and the order of the underlying group $\Z_n$ is equal to $n=105$, which is factorized as $105=3\times 5 \times 7$. Moreover, it has been proved that for any $c\in\N$, there exists a minimal vanishing sum roots of unity of order $n=3pq$ with the highest coefficient equal to $c\in\N$ for some distinct primes $p$ and $q$ (Theorem~1.1, p.~690, \citet{steinberger2008minimal}). In contrast with previous classifications based on the cardinality of the resulting sets, we choose to enumerate perfectly balanced sets according to a given order. For example, Table~\ref{tab:Expressions_N30} provides all the six expressions of the non-periodic minimal perfectly balanced sets of $\Z_{30}$ using the same notations previously introduced (p.~114 \citet{mann1965linear}), i.e. $\rho$ of order~5, $\eta$ of order 3 and $\delta$ of order 2, by also explicitly adding the associated subsets $P_1, P_2, Q_1$ and $Q_2$ related to Theorem~\ref{thm:cas2pq}.

\begin{table}[h]
\centering
\scriptsize
\caption{Expressions of the 6 non-periodic minimal perfectly balanced sets of~$\Z_{30}$ and their associated subsets $P_1$ and $P_2$ with $Q_1=1$, $Q_2=\eta + \eta^2$, $\rho^5=1$, $\eta^3=1$ and $\delta^2=1$.} 
\begin{tabular}{ccc}
\toprule
Perfectly balanced sets of $\Z_{30}$ & $P_1$ & $P_2$ \\
\midrule
$\rho + \rho^2 + \rho^3 + \rho^4 + \delta \eta + \delta \eta^2$ & $\rho + \rho^2 + \rho^3 + \rho^4$ & 1  \\
$1 + \rho^2 + \rho^3 + \delta \eta (\rho + \rho^4) + \delta \eta^2 (\rho + \rho^4)  $ & $1 + \rho^2 + \rho^3$ & $ \rho + \rho^4$  \\
$1 + \rho + \rho^4 + \delta \eta (\rho^2 + \rho^3) + \delta \eta^2 (\rho^2 + \rho^3)  $ & $1 + \rho + \rho^4$ & $ \rho^2 + \rho^3$ \\
$\rho + \rho^4 + \delta \eta (1 + \rho^2 + \rho^3) + \delta \eta^2 (1 + \rho^2 + \rho^3)  $ & $\rho + \rho^4$ & $1 +  \rho^2 + \rho^3$ \\
$\rho^2 + \rho^3 + \delta \eta (1 +  \rho + \rho^4) + \delta \eta^2 (1 +  \rho + \rho^4)$ & $\rho^2 + \rho^3$ & $1 +  \rho + \rho^4$ \\
$1  + \delta \eta (\rho + \rho^2 + \rho^3 + \rho^4) + \delta \eta^2 (\rho + \rho^2 + \rho^3 + \rho^4)$ & 1 & $\rho + \rho^2 + \rho^3 + \rho^4$ \\
\bottomrule
\end{tabular}
\label{tab:Expressions_N30}
\end{table}


\subsection{Decomposition of Minimal Perfectly Balanced Sets} \label{subsec:Algebraic_Decomposition}

As we have seen previously with Theorem~\ref{thm:redei-debuijn-schoenberg}, all non-periodic minimal perfectly balanced sets can be interpreted as an integer linear combination of regular polygons. In the particular case of the present paper, i.e. $n=2pq$, this means that they can be expressed as an integer linear combination of $p$-gons, $q$-gons and bipoints. However, the decomposition is not unique, and to our knowledge, the question of establishing all such possible decompositions remains open. The following proposition provides a list of six different decompositions for any element of $\Gamma_{2pq}$ using the previous notations.

\begin{proposition} \label{prop:decomposition}
Let $A \in \Gamma_{2pq}$. Then $A$ admits the following six different decompositions into sum and difference of subgroups of prime order of $\Z_{2pq}$:
\begin{enumerate}[label=(\roman*)]
    \item $ A = Q_1  C_p + \delta P_2  C_q - P_2 Q_1  C_2, $
    \item $ A = Q_1  C_p - P_2  C_q + P_2 Q_2  C_2, $
    \item $ A = - \delta Q_1  C_p + \delta P_2  C_q + P_1 Q_1  C_2, $
    \item $ A = \delta Q_2  C_p + P_1  C_q - P_1 Q_2  C_2, $
    \item $ A = \delta Q_2  C_p - \delta P_1  C_q + P_1 Q_1  C_2, $
    \item $ A = - Q_2  C_p + P_1  C_q + P_2 Q_2  C_2. $
\end{enumerate}
\end{proposition}

\begin{proof}
More generally, we set $n=p_1p_2p_3$ (where one of the three prime factors is equal to~2) and $(X_1,X_2)$, $(Y_1,Y_2)$ and $(Z_1,Z_2)$ as nontrivial partitions of $C_{p_1}$, $C_{p_2}$ and $C_{p_3}$. Using the notation of Remark~\ref{rem:Notation_A_P1Q1}, we have $A = X_1 Y_1 Z_1 + X_2 Y_2 Z_2$, which leads to:
$$\begin{aligned}
    A &= X_1 Y_1 Z_1 + X_2 Y_2 Z_2  \\
    &= X_1 Y_1 Z_1 + X_2 Y_1 Z_1 - X_2 Y_1 Z_1 + X_2 Y_2 Z_2 \\
    &= Y_1 Z_1 C_{p_1} - X_2 Y_1 Z_1 + X_2 Y_2 Z_2 \\
    &= Y_1 Z_1 C_{p_1} - X_2 Y_1 Z_1 + X_2 Y_2 Z_2 + X_2 Y_1 Z_2 - X_2 Y_1 Z_2  \\
    &= Y_1 Z_1 C_{p_1} + X_2 Z_2 C_{p_2} - X_2 Y_1 Z_1 - X_2 Y_1 Z_2  \\
    &= Y_1 Z_1 C_{p_1} + X_2 Z_2 C_{p_2} - X_2 Y_1 C_{p_3}
\end{aligned}$$
From this last expression of $A$, it follows that each decomposition is completely determined by the ordering of the three prime numbers 2, $p$ and $q$ assigned to $p_1$, $p_2$ and $p_3$. We conclude that there are six possible expressions corresponding to the number of permutations.
\end{proof}

Note that it is possible to express all these six decompositions into a unique formula which highlights the role of the subsets underlying the decomposition. Using the notations defined in the previous proof, any non-periodic minimal perfectly balanced set $A$ can be expressed as follows:
$$ A = Y_1 Z_1  C_{p_1} + X_2 Z_2  C_{p_2} - X_2 Y_1  C_{p_3}.$$
For example, in the decomposition $(i)$ of Proposition~\ref{prop:decomposition}, we have $(p_1,p_2,p_3)=(p,q,2)$ and $(Z_1,Z_2)=(1,\delta)$, but in the decomposition $(ii)$ of Proposition~\ref{prop:decomposition} we have $(p_1,p_2,p_3)=(p,2,q)$ and $(Y_1,Y_2)=(1,\delta)$. This general formula has the advantage of compacting the six different decompositions of Proposition~\ref{prop:decomposition} into a unique expression. Moreover, it is possible to obtain the value of the lowest and highest cardinalities of all the non-periodic minimal perfectly balanced sets of a given order $2pq$. They are given by the following proposition.

\begin{proposition} \label{prop:highest_lowest_cardinalities}
The value of the lowest and highest cardinalities of the elements of $\Gamma_{2pq}$ are respectively equal to: $$ p+q-2 \; (\text{lowest}) \quad \text{and} \quad pq -p -q +2 \; (\text{highest}).$$
\end{proposition} 

\begin{proof}
Let $A\in\Gamma_{2pq}$. Using the decomposition $(i)$ of Proposition~\ref{prop:decomposition}, we have:
$$|A|=|Q_1|p+|P_2|q-2|Q_1||P_2|,$$
where $|Q_1|$ (respectively $|P_2|$) takes values ranging from 1 to $q-1$ (respectively from $1$ to $p-1$). Since this expression is affine in each variable separately, its maximum and minimum are obtained at the endpoints of these intervals. For $|Q_1| = 1$ and $|P_2| = 1$ (as well as $|Q_1| = q - 1$ and $|P_2| = p - 1$), we have $|A| = p + q - 2$. On the other hand, for $|Q_1| = q - 1$ and $|P_2| = 1$ (as well as $|Q_1| = 1$ and $|P_2| = p - 1$), we have $|A| = pq - p - q + 2$. Therefore, there are only two possible extremal values. To determine which of the two is the greatest, we have: $(pq-p-q+2) - (p+q-2) = (p-2)(q-2)>0$ because $p$ and $q$ are prime greater than 2. This shows that $|A|$ achieves a maximum value equal to $pq-p-q+2$ and a minimum value equal to $p+q-2$.
\end{proof}

\begin{figure}[h!]
    \centering
    \begin{subfigure}{\textwidth}
        \centering
        \includegraphics[width=\linewidth]{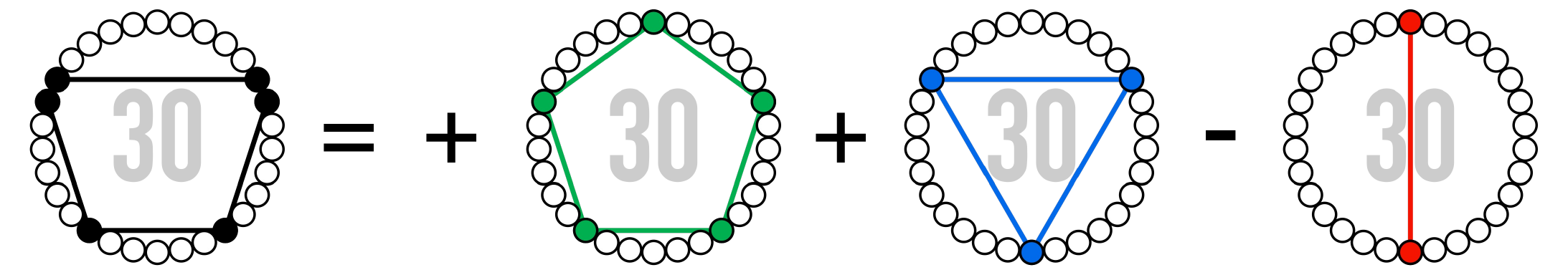}
        \caption{$A = C_5 + z^{15}  C_3 - C_2$}
        \label{subfig:Decomposition_N30_1}
    \end{subfigure} \par\vspace{0.28cm}
    \begin{subfigure}{\textwidth}
        \centering
        \includegraphics[width=\linewidth]{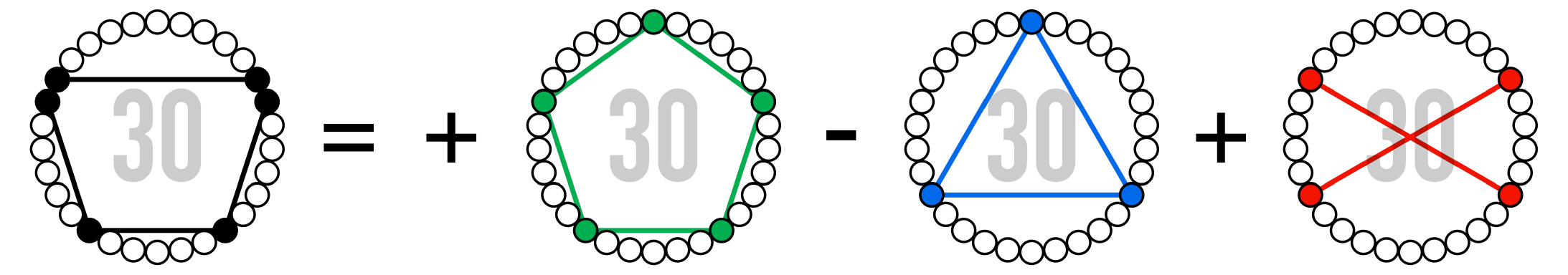}
        \caption{$A = C_5 - C_3 + (z^{10}+z^{20}) C_2$}
        \label{subfig:Decomposition_N30_2}
    \end{subfigure} \par\vspace{0.28cm}
    \begin{subfigure}{\textwidth}
        \centering
        \includegraphics[width=\linewidth]{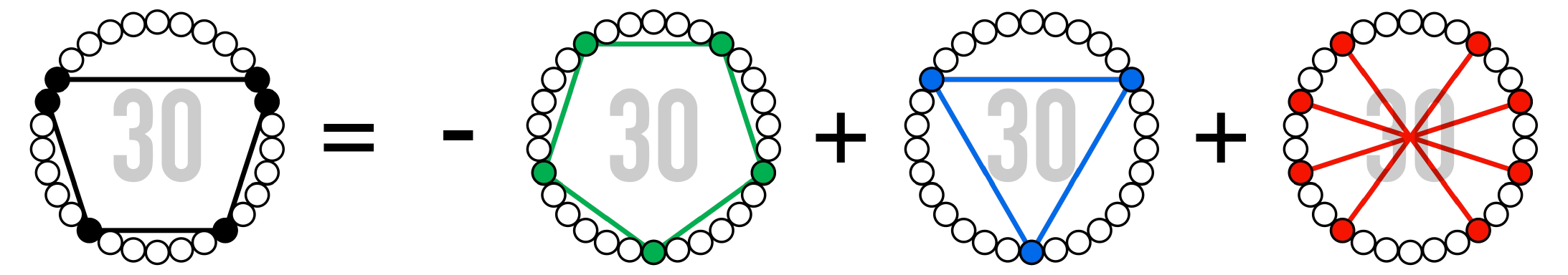}
        \caption{$A = - z^{15} C_5 + z^{15} C_3 + (z^6+z^{12}+z^{18}+z^{24}) C_2$}
        \label{subfig:Decomposition_N30_3}
    \end{subfigure} \par\vspace{0.28cm}
    \begin{subfigure}{\textwidth}
        \centering
        \includegraphics[width=\linewidth]{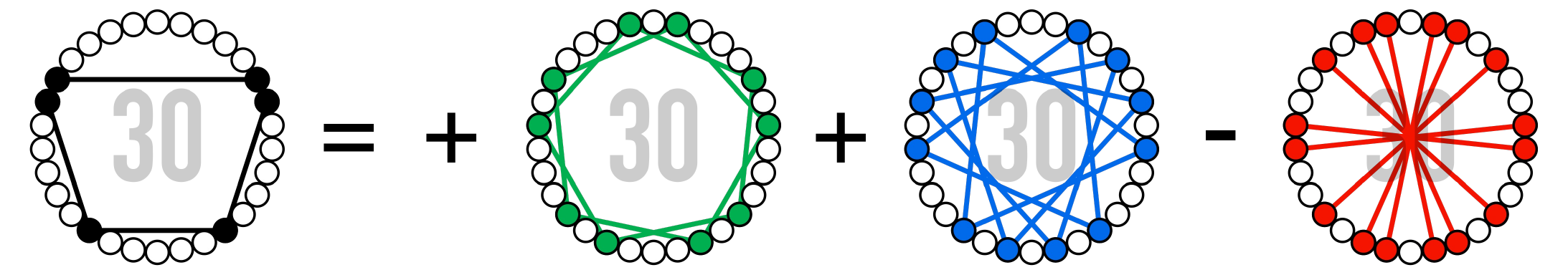}
        \caption{$A = (z^{5}+z^{25})  C_5 + (z^6+z^{12}+z^{18}+z^{24})  C_3 - (z^2+z^4+z^8+z^{14}+z^{16}+z^{22}+z^{26}+z^{28})  C_2$}
        \label{subfig:Decomposition_N30_4}
    \end{subfigure} \par\vspace{0.28cm}
    \begin{subfigure}{\textwidth}
        \centering
        \includegraphics[width=\linewidth]{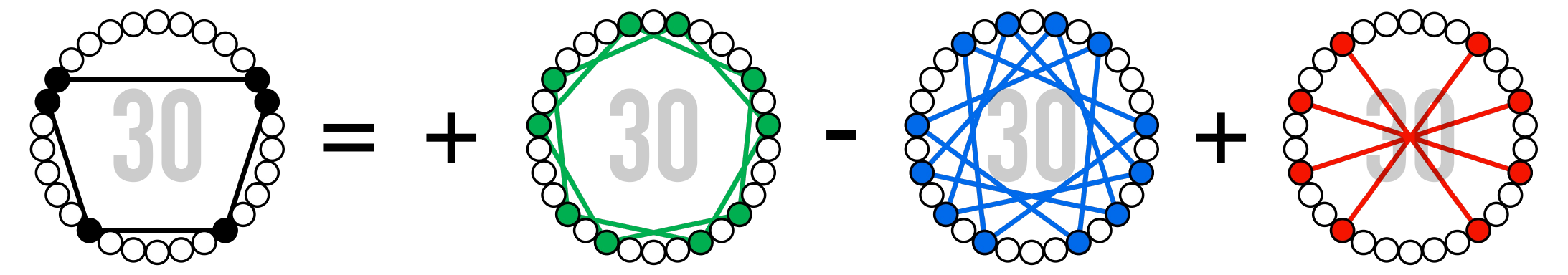}
        \caption{$A = (z^{5}+z^{25})  C_5 - (z^3+z^9+z^{21}+z^{27})  C_3 + (z^6+z^{12}+z^{18}+z^{24})  C_2$}
        \label{subfig:Decomposition_N30_5}
    \end{subfigure} \par\vspace{0.28cm}
    \begin{subfigure}{\textwidth}
        \centering
        \includegraphics[width=\linewidth]{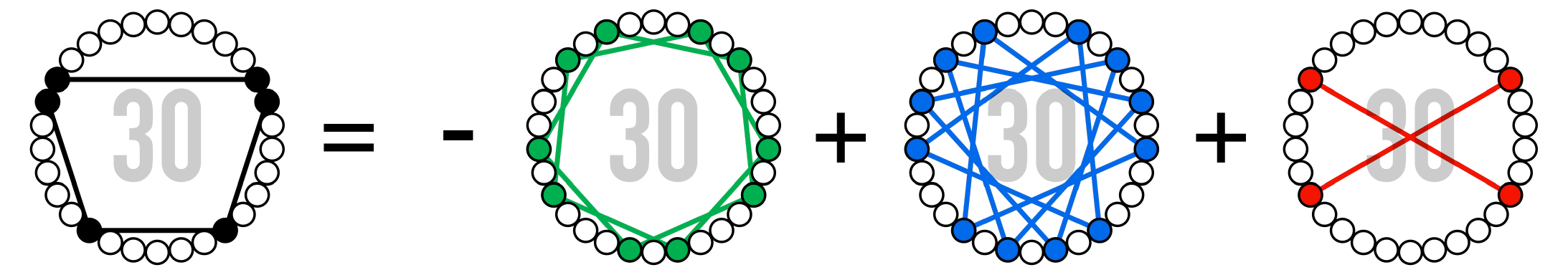}
        \caption{$A = -(z^{10}+z^{20})  C_5 + (z^6+z^{12}+z^{18}+z^{24})  C_3 + (z^{10}+z^{20})  C_2$}
        \label{subfig:Decomposition_N30_6}
    \end{subfigure}
    \caption{Visualization of six decompositions of the non-periodic minimal perfectly balanced set $A = z^5+z^6+z^{12}+z^{18}+z^{24}+z^{25}$ of $\Z_{30}$ as an integer linear combination of regular pentagons, triangles and bipoints.}
    \label{fig:Decompositions}
\end{figure}

According to Proposition~\ref{prop:highest_lowest_cardinalities}, among all the non-periodic minimal perfectly balanced sets of $\Z_{30}$ illustrated in Figure~\ref{fig:N30}, we can verify that the lowest and highest values of cardinality are respectively equal to $3+5-2=6$ (Figure~\ref{subfig:N30_1}) and $3\times5-3-5+2=9$ (Figure~\ref{subfig:N30_4}). Moreover, by considering the generator $z$ of $\Z_{30}$, the non-periodic minimal perfectly balanced set with cardinality equal to 6 (Figure~\ref{subfig:N30_1}) can be expressed as $A = z^5+z^6+z^{12}+z^{18}+z^{24}+z^{25}$. Figure~\ref{fig:Decompositions} illustrates the six different decompositions of Proposition~\ref{prop:decomposition} with $P_1= z^6+z^{12}+z^{18}+z^{24}$, $P_2 = 1$, $Q_1 = 1$ and $Q_2=z^{10}+z^{20}$. In the case of $\Z_{30}$, they correspond to integer linear combinations of regular pentagons $C_5$, triangles $C_3$ and bipoints $C_2$.

\begin{remark}
It is important to observe that in each of the six decompositions of Proposition~\ref{prop:decomposition}, all coefficients appearing in front of the terms $C_p$, $C_q$, and $C_2$ belong to $\{-1,0,1\}$ and they have the same sign with respect to $C_p$, $C_q$, and $C_2$. Indeed, each coefficient arises from a product $X_iY_j$, where $X_i$ and $Y_j$ are subsets of subgroups of coprime orders (for instance $P_2Q_1$ with $P_2$ and $Q_1$ as subset of order $p$ and $q$). Such a product is itself a set, i.e. with coefficients in $\{0,1\}$. Thus, all coefficients with respect to each term $C_p$, $C_q$, and $C_2$ belong to $\{0,1\}$, or $\{-1,0\}$ when there is a minus sign, and therefore the sign remains constant with respect to $C_p$, $C_q$, and $C_2$. For example, the decomposition of $A$ in Figure~\ref{subfig:Decomposition_N30_6} is $ A = - (z^{10}+z^{20})  C_5 + (z^6+z^{12}+z^{18}+z^{24})  C_3 + (z^{10}+z^{20})  C_2 $ and the coefficients in front of $C_5$ are all equal to 0 or $-1$, while the coefficients in front of $C_3$ and $C_2$ are all equal to 0 or 1. This property confirms the behavior conjectured by A. Milne et al. according to the coefficients in front of the prime polygon patterns (p.~114, \citet{milne2017exploring}) that we prove here with these decompositions in the particular case of order $2pq$.
\end{remark}


\section{Combinatorial Enumeration of Perfectly Balanced Sets} \label{sec:Combinatorial}

In this section, we investigate the number and the distribution of non-periodic minimal perfectly balanced sets of order $2pq$. We first derive a formula for counting these sets. We then show that their distribution is symmetric with respect to the cardinality of the sets. Finally, we propose a formula for the possible cardinalities and the distribution of non-periodic minimal perfectly balanced sets in the particular case of order $6p$ (i.e. when $q=3$).


\subsection{Number of Minimal Perfectly Balanced Sets} \label{subsec:Combinatorial_Number}

We propose here to count the number of elements of $\Gamma_{2pq}$ up to rotation. It is known that there are exactly three periodic minimal perfectly balanced sets up to rotation, respectively $C_2$, $C_p$, and $C_q$. The following formula gives the number of non-periodic ones.

\begin{proposition} \label{prop:Number_of_PBS}
The number of elements of $\Gamma_{2pq}$ up to rotation is: $$\frac{2\,(2^{p-1} - 1)(2^{q-1} - 1)}{pq}.$$
\end{proposition}

\begin{proof}
Using Theorem~\ref{thm:cas2pq}, it suffices to count the possible choices of $P_1$ and $Q_1$. Since $|C_p| = p$ and $|C_q| = q$, there are $2^p$ subsets of $C_p$ and $2^q$ subsets of $C_q$. Excluding the empty sets and the full sets, we obtain $2^p - 2$ choices for $P_1$ and $2^q - 2$ choices for $Q_1$. Therefore, there are
$4 (2^{p-1} - 1)(2^{q-1} - 1)$ non-periodic minimal perfectly balanced sets of order $2pq$. To avoid counting rhythms that are equivalent up to rotation, we divide by $2pq$ which is the number of possible rotations because the sets are non-periodic.
\end{proof}

Table~\ref{tab:enumeration_results} gives the number of non-periodic minimal perfectly balanced sets of order $2pq$ up to rotation, along with the corresponding order, for the first few values of $p$ and $q$. Note that whereas determining the total number of sets is immediate, the actual enumeration of these sets by brute force methods can be cumbersome and slow given the need to compute all possible subsets of $C_p$ and $C_q$.

\begin{table}[h]
\centering
\scriptsize
\captionsetup{margin=20pt}
\caption{Number of non-periodic minimal perfectly balanced sets of order $2pq$ up to rotation for the first values of $p$ and~$q$. The value $2pq$ is indicated in parentheses.} 
\begin{tabular}{ccccccc}
\toprule
\diagbox[width=3em]{$q$}{$p$} & $3$ & $5$ & $7$ & $11$ & $13$ & $17$ \\
\midrule
$3$  & -- & 6(30) & 18(42) & 186(66) & 630(78) & 7710(102) \\
$5$  &  --  & -- & 54(70) & 558(110) & 1890(130) & 23130(170) \\
$7$  &  --  &  --  & -- & 1674(154) & 5670(182) & 69390(238) \\
$11$ &  --  & --   &  --  & -- & 58590(286) & 717030(374) \\
\bottomrule
\end{tabular}
\label{tab:enumeration_results}
\end{table}


\subsection{Distribution of Minimal Perfectly Balanced Sets} \label{subsec:Combinatorial_Distribution}

We now investigate the distribution of non-periodic minimal perfectly balanced sets by classifying them according to their cardinality and reveal that this distribution is symmetric when $n=2pq$. We have already presented the 6 non-periodic minimal perfectly balanced sets of $\Z_{30}$ in Figure~\ref{fig:N30}. Their distribution has a clear symmetric character: there is 1 set of cardinality $6$, 2 sets of cardinality $7$, 2 sets of cardinality $8$, and 1 set of cardinality~$9$. The same phenomenon holds for the $18$ non-periodic minimal perfectly balanced sets of $\Z_{42}$ shown in Figure~\ref{fig:N42} in the introduction section. In particular, there are respectively $1$, $3$, $5$, $5$, $3$, and $1$ sets with cardinalities ranging from $8$ to $13$. In this case, it is worth noting that two of the sets of cardinalities $10$ and $11$, which are displayed on the right of Figure~\ref{fig:N42}, are not \textit{symmetric}. This can be seen from the fact that they have no axis of symmetry in the circular representation, or more formally, based on the following definition.

\begin{definition}
Let $A \in \Z[C_n]$. The element $A=\sum_{i=0}^{n-1}a_iz^i$ is \textit{symmetric} if there exists $t\in C_n$ such that $A=t\bar{A}$, where $\bar{A}=\sum_{i=0}^{n-1}a_iz^{-i}$.
\end{definition}

The possibility of having non-periodic minimal perfectly balanced sets that are not symmetric for $\Z_{42}$ contrasts with the case of $\Z_{30}$, where all non-periodic minimal perfectly balanced sets are symmetric. The following proposition shows that $\Z_{30}$ is the only case in which this property holds.

\begin{proposition}
All non-periodic minimal perfectly balanced sets of order $n=2pq$ are symmetric iff $n=30$.
\end{proposition}

\begin{proof}
If $n=30$, it is clear that all the sets are symmetric (see Figure~\ref{fig:N30}). For $n > 30$, it suffices to show that there always exists a non-periodic minimal perfectly balanced set $A$ which is non-symmetric. Without loss of generality, suppose that $2<p<q$, which implies that $q\geq7$. We set $P_1=1$ and $Q_1=1+\eta+\eta^3$, which leads to:
$$ A = 1+\eta+\eta^3 + \delta P_2 Q_2 \quad \text{ and } \quad  \bar{A} = 1+\eta^{q-1}+\eta^{q-3} + \delta \overline{P_2 Q_2}.$$

Assuming that $A$ is symmetric, there exists $t\in C_n$ such that $A=t\bar{A}$. The cardinalities of $1 + \eta + \eta^3$ and $1+\eta^{q-1}+\eta^{q-3}$ are 3 and the cardinalities of $\delta P_2 Q_2$ and $\delta \overline{P_2 Q_2}$ are $|P_2|\times |Q_2|=(p-1)(q-3)>3$. Therefore, because $\Z_{2pq} \simeq \Z_{2} \times \Z_{pq}$ and by a cardinality argument, we have: 
$$ 1+\eta+\eta^3 = t(1+\eta^{q-1}+\eta^{q-3}),$$ 
which is absurd because $1+\eta+\eta^{3}$ is not equivalent to $1+\eta^{q-1}+\eta^{q-3}$ up to rotation since $q\geq7$. This proves that $A$ is not symmetric. Therefore, $n=30$ is the only subdivision for which all non-periodic minimal perfectly balanced sets are symmetric.
\end{proof}

Before demonstrating the symmetric character of the distribution of the non-periodic minimal perfectly balanced sets of order $2pq$, we now focus on the special case of order $6p$, i.e. when $q=3$. In particular, the following proposition provides the cardinality of the sets and gives their number according to a specific cardinality. It reveals that the distribution follows the binomial coefficients and the cardinality takes $p-1$ values, ranging from $p+1$ to $2p-1$. Similar results have been obtained concerning only the lower bound of the number of elements of a vanishing sum of roots of unity (Theorem 5.2, p.~103, \citet{lam2000vanishing}) with recent improvements in a particular case (Theorem 1.2, p.~3, \citet{Laba2022Vanishing}). We provide here a more specific description of all the possible cardinalities of non-periodic minimal perfectly balanced sets of order $6p$.

\begin{proposition} \label{proposition:distribution_6p}
Up to rotation, the number of elements of $\Gamma_{6p}$ with cardinality $p+j$ (i.e. from $p+1$ to $2p-1$) is given by:
$$ \frac{\binom{p}{j}}{p}, \qquad \forall j = 1, \dots, p-1. $$
\end{proposition} 

\begin{proof}
Starting from the following description of $A \in \Gamma_{2pq}$ as $P_1 Q_1 + \delta P_2 Q_2$, we derive its cardinality: $$ |A| = |P_1| \times |Q_1| + |P_2| \times |Q_2|.$$
In the particular case of $\Gamma_{6p}$, we consider $(Q_1,Q_2)$ a nontrivial partition of $C_3$ where $|Q_1| = 1 $ and $ |Q_2|=2 $ up to a permutation of $Q_1$ and $Q_2$. Therefore, we have: $$|A| = |P_1| + 2 \times |P_2|,$$
which shows that the cardinality of $A$ goes from $p+1$ (when $|P_1|=p-1$ and $|P_2|=1$) to $2p-1$ (when $|P_1|=1$ and $|P_2|=p-1$) with an incremental step equal to 1. Therefore, by supposing $|Q_1| = 1 $ and $ |Q_2|=2 $, the only possibility of having $|A|=p+k$ is $|P_2|=k$. There are $\binom{p}{k}$ choices for $P_2$ with $k$ elements and $3$ choices for $Q_1$ with 1 element. By a similar reasoning on $P_1$ and $Q_2$, there are $\binom{p}{k}\times 3 \times 2$ different sets $A$ of cardinality $p+k$. Since these sets are all non-periodic, there are exactly $6p$ rotations, which leads to $\binom{p}{k}\times 3 \times 2/6p=\binom{p}{k}/p$ non-periodic minimal perfectly balanced sets of cardinality $p+k$ of order $6p$ up to rotation.
\end{proof}

\vspace{0.3cm}

Table~\ref{table:distribution} provides the distribution of the elements of $\Gamma_{2pq}$ up to rotation. Note that some values, such as $70=2\times 5 \times 7$, are not of the form $6p$. Therefore, despite its symmetric character, the distribution does not respect the behaviour described in Proposition~\ref{proposition:distribution_6p}: the incremental step of the cardinalities is not equal to 1 (for example it goes from 10 to 13), and the number of sets of the first half is not increasing (for example it goes from 3 to 2). However, in all cases, for a given order $2pq$, the cardinality of these sets are bounded by the lowest and highest values, respectively $p+q-2$ and $pq-p-q+2$, described in Proposition~\ref{prop:highest_lowest_cardinalities}.

\begin{table}[h]
\caption{Number of non-periodic minimal perfectly balanced sets up to rotation ranked by cardinality (in parentheses) for different values of order $2pq$. The values 30, 42, 66, 78, 102, 114 and 138 are of the form $6p$ ($q=3$) and therefore follow the distribution described in Proposition~\ref{proposition:distribution_6p}, while 70, 110 and 130 are of the form $10p$ ($q=5$) and 154 is of the form $14p$ ($q=7$).}
\centering
\scriptsize
\begin{tabular}{ccccccccccc}
\toprule
$30$ & $42$ & $66$ & $70$ & $78$ & $102$ & $110$ & $114$ & $130$ & $138$ & $154$ \\
\midrule
1(6) & 1(8) & 1(12) & 1(10) & 1(14) & 1(18) & 1(14) & 1(20) & 1(16) & 1(24) & 1(16) \\
2(7) & 3(9) & 5(13) & 3(13) & 6(15) & 8(19) & 5(17) & 9(21) & 6(19) & 11(25) & 5(21) \\
2(8) & 5(10) & 15(14) & 2(15) & 22(16) & 40(20) & 15(20) & 51(22) & 22(22) & 77(26) & 3(25) \\
1(9) & 5(11) & 30(15) & 11(16) & 55(17) & 140(21) & 32(23) & 204(23) & 55(25) & 385(27) & 15(26) \\
      & 3(12) & 42(16) & 10(17) & 99(18) & 364(22) & 10(24) & 612(24) & 2(27) & 1463(28) & 15(28) \\
      & 1(13) & 42(17) & 10(18) & 132(19) & 728(23) & 30(25) & 1428(25) & 111(28) & 4389(29) & 75(31) \\
      &       & 30(18) & 11(19) & 132(20) & 1144(24) & 102(26) & 2652(26) & 44(29) & 10659(30) & 95(34) \\
      &       & 15(19) & 2(20) & 99(21) & 1430(25) & 84(27) & 3978(27) & 110(30) & 21318(31) & 25(35) \\
      &       & 5(20) & 3(22) & 55(22) & 1430(26) & 84(28) & 4862(28) & 330(31) & 35530(32) & 117(36) \\
      &       & 1(21) & 1(25) & 22(23) & 1144(27) & 102(29) & 4862(29) & 264(32) & 49742(33) & 276(37) \\
      &       &       &       & 6(24) & 728(28) & 30(30) & 3978(30) & 264(33) & 58786(34) & 210(38) \\
      &       &       &       & 1(25) & 364(29) & 10(31) & 2652(31) & 330(34) & 58786(35) & 210(39) \\
      &       &       &       &       & 140(30) & 32(32) & 1428(32) & 110(35) & 49742(36) & 276(40) \\
      &       &       &       &       & 40(31) & 15(35) & 612(33) & 44(36) & 35530(37) & 117(41) \\
      &       &       &       &       & 8(32) & 5(38) & 204(34) & 111(37) & 21318(38) & 25(42) \\
      &       &       &       &       & 1(33) & 1(41) & 51(35) & 2(38) & 10659(39) & 95(43) \\
      &       &       &       &       &       &       & 9(36) & 55(40) & 4389(40) & 75(46) \\
      &       &       &       &       &       &       & 1(37) & 22(43) & 1463(41) & 15(49) \\
      &       &       &       &       &       &       &  & 6(46) & 385(42) & 15(51) \\
      &       &       &       &       &       &       &  & 1(49) & 77(43) & 3(52) \\
      &       &       &       &       &       &       &  &       & 11(44) & 5(56) \\
      &       &       &       &       &       &       &  &       & 1(45) & 1(61) \\
\bottomrule
\end{tabular}
\label{table:distribution}
\end{table}


\newpage

\section{Duality of Perfectly Balanced Sets} \label{sec:Pairwise}

In this section, we present a major property which explains the symmetric distribution of perfectly balanced sets. We called this property the \textit{duality} of non-periodic minimal perfectly balanced sets of order~$2pq$. It reveals that any element of $\Gamma_{2pq}$ can be combined with another one to produce a self-complementary $2p$-periodic set. Moreover, we prove that the reciprocal also holds, i.e. any self-complementary $2p$-periodic set can be decomposed into two elements of $\Gamma_{2pq}$ that we call \textit{dual sets}. To our knowledge, this property was not noticed in the literature. In addition, we further characterize the duality in the special case of order $6p$ by describing the two possible configurations of the resulting self-complementary set: if it is $6$-periodic, it is unique up to rotation, whereas if it is $2p$-periodic, it admits a unique decomposition into dual non-periodic minimal perfectly balanced sets.


\subsection{Main Results on the Duality of Minimal Perfectly Balanced Sets} \label{subsec:Pairwise_Main}

Before presenting the main results, we introduce the concept of \textit{self-complementary} in the case of sets.

\begin{definition}
Let $A \in \Z[C_{2pq}]$ be a set. Then, $A$ is \textit{self-complementary} if $A^c=\delta A$, where $A^c=C_n-A$ is the \textit{complement} of $A$ and $\delta$ is the element of order 2 of $C_{2pq}$.
\end{definition}

We now present the main result of this paper which reveals that non-periodic minimal perfectly balanced sets of order~$2pq$ are always related by duality, although two dual sets do not have the same cardinality. Their association is similar to a complementary relation in the sense that it produces a set of $\Z_{2pq}$ which is $2p$-periodic and self-complementary.

\begin{theorem} \label{thm:Pairwise_Association}
Let $A\in\Gamma_{2pq}$. Then, there exists $B\in\Gamma_{2pq}$ such that: 
$$A + B = C,$$
where $C$ is a subset of $\Z_{2pq}$ which is $2p$-periodic and self-complementary.
\end{theorem}

\begin{proof}
Since $A$ is a non-periodic minimal perfectly balanced set of $\Z_{2pq}$, we can use Theorem~\ref{thm:cas2pq}, which allows to describe $A$ as follows:
$$A=P_1Q_1+\delta P_2Q_2,$$
with the nontrivial partitions $(P_1,P_2)$ and $(Q_1,Q_2)$ of $C_p$ and $C_q$. We then define:
$$ B = P_1Q_2+\delta P_2Q_1 .$$
The set $B$ is indeed a non-periodic minimal perfectly balanced set of order $2pq$ by Theorem~\ref{thm:cas2pq}. Setting $A+B=C$, we therefore have: $$C= C_q (P_1 + \delta P_2).$$

We first show that $C$ is a subset of $\Z_{2pq}$. Because $C_qP_1$ and $C_q\delta P_2$ are both a product of sets of coprime order, they are sets. Moreover, since $P_1\cap\delta P_2=\emptyset$, the sum of $C_qP_1$ and $C_q\delta P_2$ has coefficients equal to $0$ or $1$. Hence, $C$ is a set.

Next, we show that $C$ is a self-complementary set. We have: 
$$ C + \delta C = C_q(P_1+\delta P_2)+ C_q(\delta P_1+P_2) = C_qC_pC_2 = C_{2pq},$$
which shows that $C^c = \delta C$ and therefore $C$ is self-complementary.

Finally, we show that $C$ is $2p$-periodic. It is clear that $C$ is periodic with period at least $2pq/q = 2p$, since $C_q$ has this property. We then show that this is the smallest period. Since any period of $C$ must divide $2p$, it is enough to show that $C$ has not a period equal to $p$ (i.e. $\delta \eta C \neq C$ because $\delta \eta$ has order $2q$) and $2$ (i.e. $\rho \eta C \neq C$ because $\rho \eta$ has order $pq$). Since $\eta C_q = C_q$, we have $\eta C = C$, hence we only need to show that $\delta C \neq C$ and $\rho C \neq C$. The first inequality come directly from the self-complementarity of $C$. We prove the second inequality by contradiction.
$$\begin{aligned}
& & \rho C = C && (\text{assumption}) \\
& \Rightarrow & \rho C_q (P_1 + \delta P_2) = C_q (P_1 + \delta P_2) && (\text{definition of } C) \\
& \Rightarrow & C_q (\rho P_1 + \delta \rho  P_2) = C_q (P_1 + \delta P_2) && (\text{commutativity}) \\
& \Rightarrow & \rho P_1 + \delta \rho  P_2 = P_1 + \delta P_2 && (\Z_{2pq} \simeq \Z_{q} \times \Z_{2p}) \\
& \Rightarrow & \rho P_1 = P_1 \text{ and } \rho P_2 = P_2 && (\Z_{2p} \simeq \Z_{2} \times \Z_{p}) \\
& \Rightarrow & P_1 = \emptyset \text{ or } P_1 = C_p && (\rho^p=1) \\
\end{aligned}$$
which is a contradiction because $(P_1,P_2)$ is a nontrivial partition of $C_p$.
\end{proof}

Theorem~\ref{thm:Pairwise_Association} explains the symmetric character of the distribution of the elements of $\Gamma_{2pq}$ represented in Table~\ref{table:distribution}. In fact, given such a set $A$ of cardinality $k$, there exists a non-periodic minimal perfectly balanced set $B$ of cardinality $k'$ such that $k+k'=pq$, because $A+B$ is a self-complementary set of $\Z_{2pq}$ and therefore has cardinality $pq$. A direct consequence of this result is that any non-periodic minimal perfectly balanced sets of order $2pq$ has cardinality less than $pq$, which is confirmed by the exact highest cardinality equal to $pq-p-q+2$ (Proposition~\ref{prop:highest_lowest_cardinalities}).

\begin{figure}[h]
     \centering
     \includegraphics[width=0.62\textwidth]{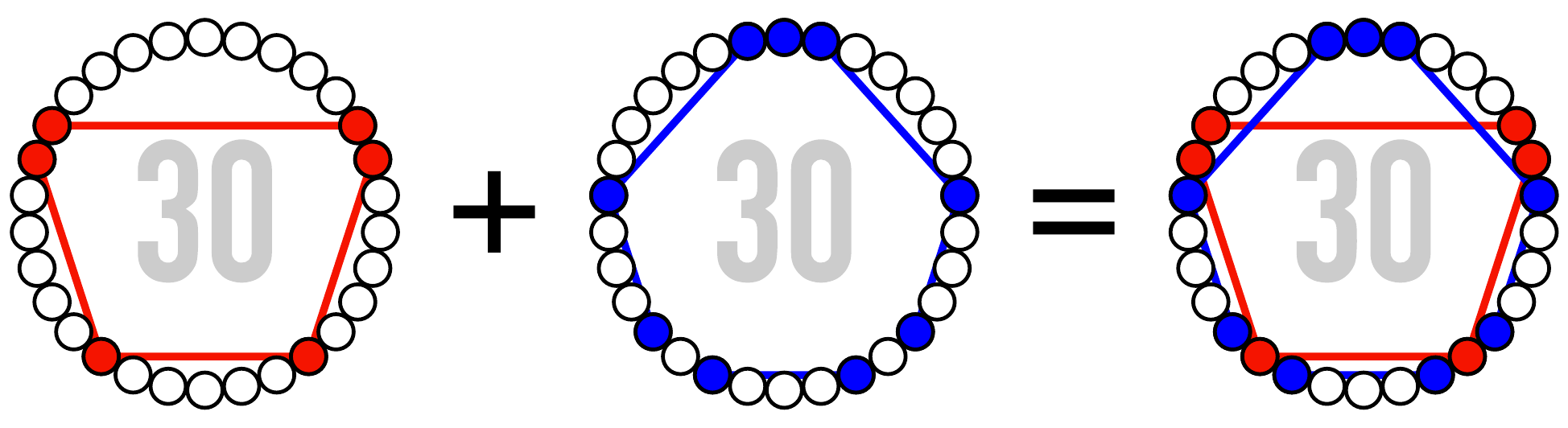} \\
     \includegraphics[width=0.62\textwidth]{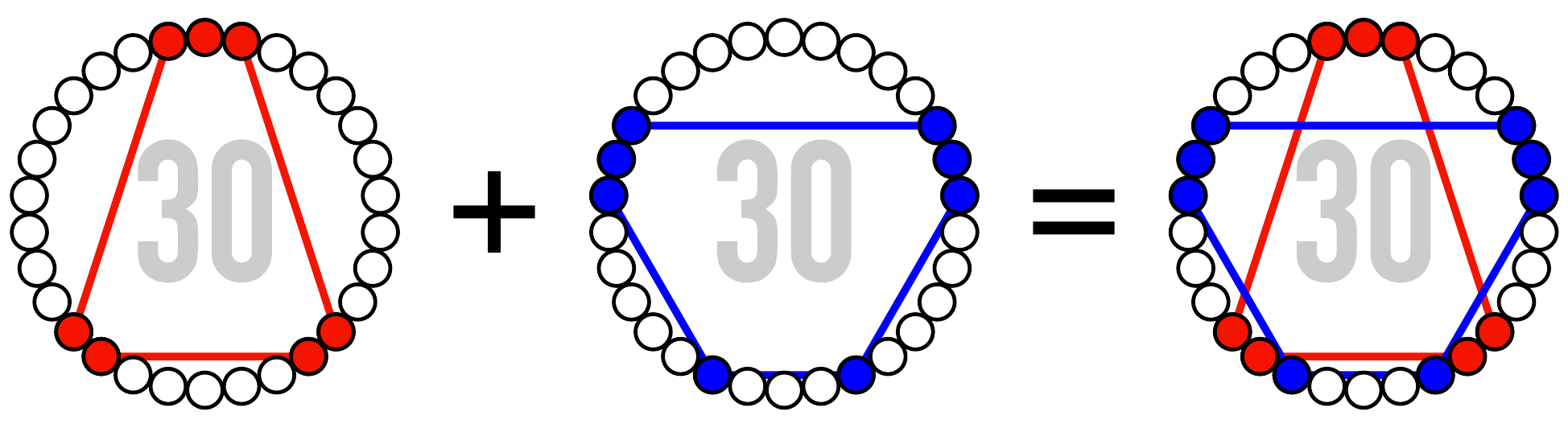} \\
     \includegraphics[width=0.62\textwidth]{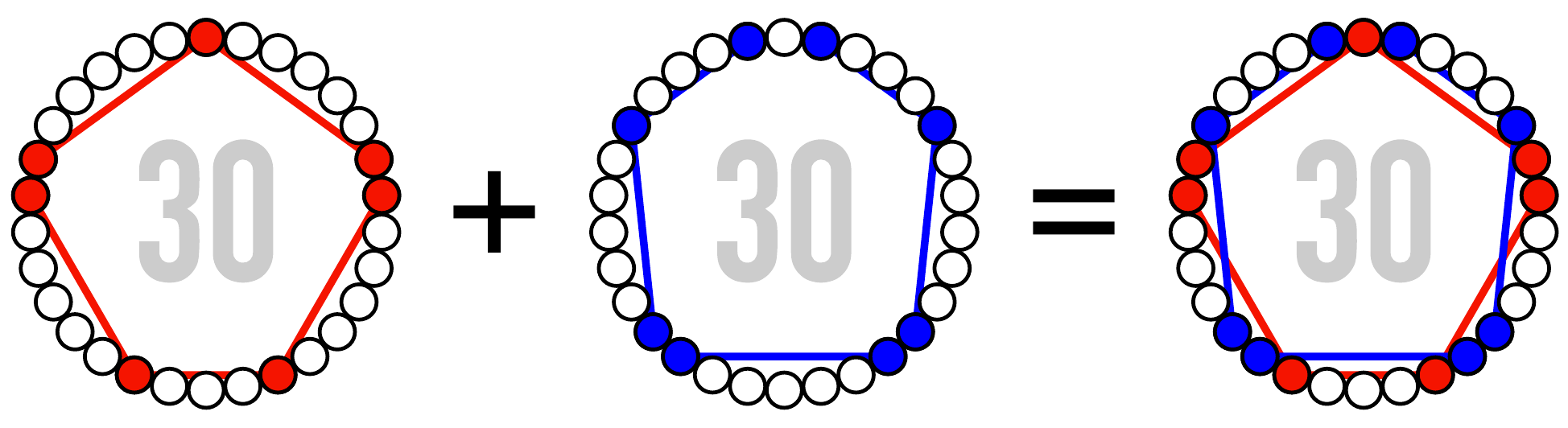}
     \caption{The dual associations of non-periodic minimal perfectly balanced sets $A+B$ of $\Z_{30}$ to form the same self-complementary $6$-periodic set.}
     \label{fig:Association_30_1}
\end{figure}

Figure~\ref{fig:Association_30_1} illustrates the three different possible dual associations of $\Z_{30}$ to obtain a self-complementarity 6-periodic set. Note that, in this case, the resulting set is always the same (corresponding of five concatenations  of the pattern $111000$).

\begin{remark}
It is possible to exchange the roles of $p$ and $q$ to obtain a dual association that is $2q$-periodic (rather than $2p$-periodic). In particular, by using the notation of Theorem~\ref{thm:Pairwise_Association}, and writing $\delta$ as the element of order 2, we also have: 
$A + \delta B = C',$
where $C'=C_p(Q_1+\delta Q_2)$ is a $2q$-periodic and self-complementary set of $\Z_{2pq}$ by the same arguments used in the proof of Theorem~\ref{thm:Pairwise_Association}.
\end{remark}

\begin{figure}[h!]
     \centering
     \includegraphics[width=0.62\textwidth]{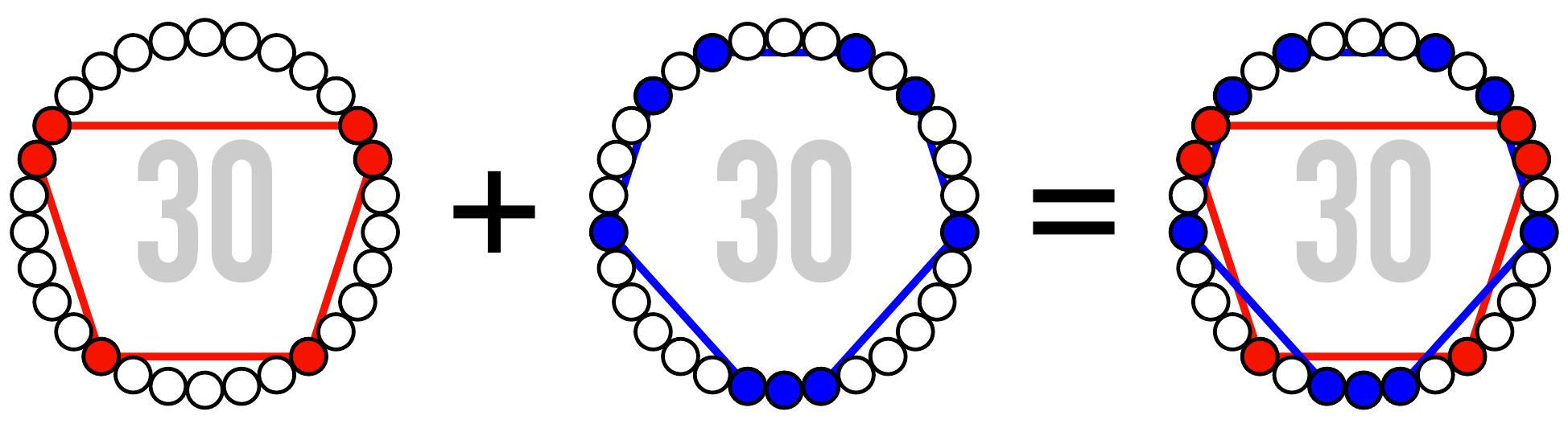} \\
     \includegraphics[width=0.62\textwidth]{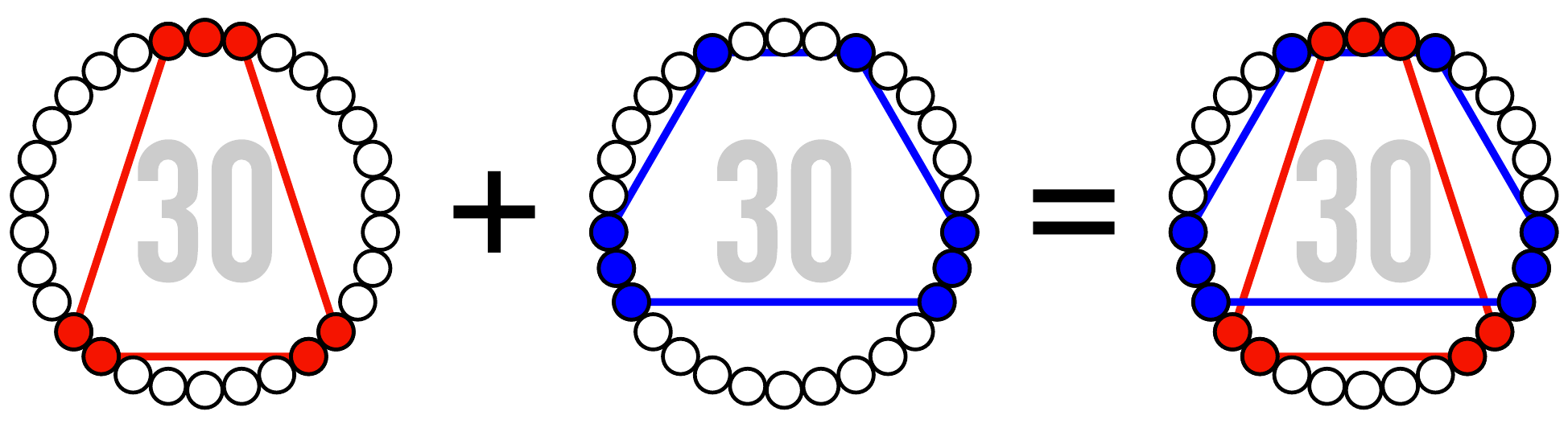} \\
     \includegraphics[width=0.62\textwidth]{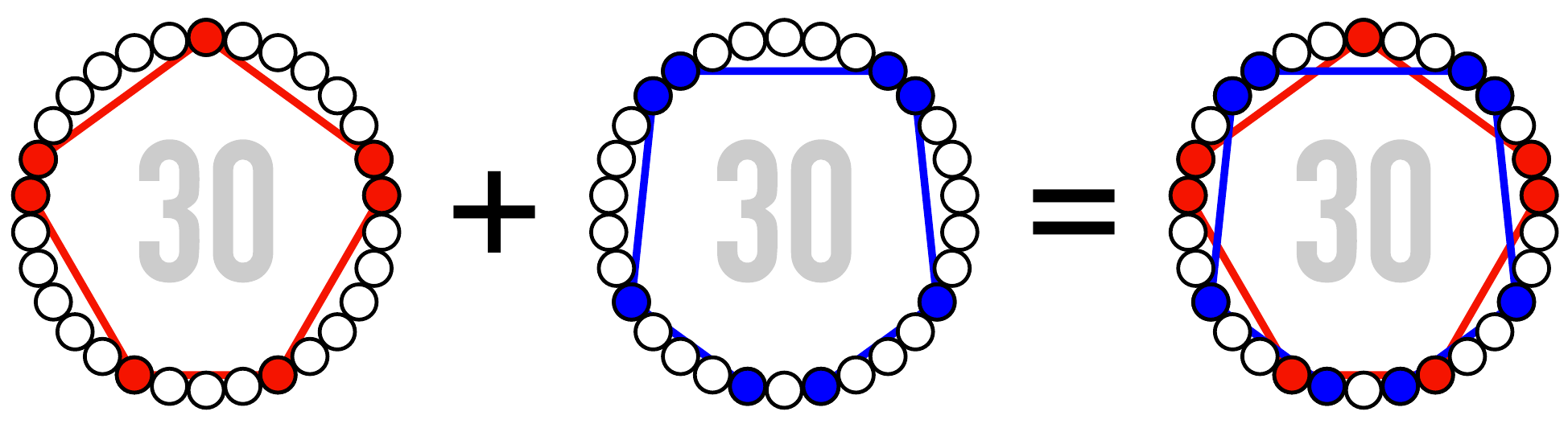}
     \caption{Three different self-complementary $10$-periodic sets of $\Z_{30}$ obtained by starting from the same perfectly balanced sets of Figure~\ref{fig:Association_30_1} but by applying a rotation by $180^\circ$ to the dual sets leading to $A+\delta B$.}
     \label{fig:Association_30_2}
\end{figure}

The previous remark implies that there also exist three dual associations of $\Z_{30}$ leading to self-complementary 10-periodic sets (instead of 6-periodic ones). These three dual associations are illustrated in Figure~\ref{fig:Association_30_2}. 

Note that the resulting sets are all different. They can be expressed as three concatenations of patterns of length 10, which are equal respectively to: 1110100010 (first association), 1111100000 (second association) and 1101100100 (third association). Surprisingly, the dual associations $A+B$ of Figure~\ref{fig:Association_30_1} and Figure~\ref{fig:Association_30_2} are related by a very simple transformation: the $B$ sets are simply rotated by $180^\circ$. It is worth noticing that these rotations produce periodic sets which are different in their inner structure and periodicities. Moreover, this property holds in the general case of order $2pq$: given any dual association, it is possible to obtain another one by simply rotating one of the two sets. This changes the resulting set and in particular its periodicity.

\subsection{Decomposition of a Self-Complementary Set Into Minimal Perfectly Balanced Sets} \label{subsec:Reciprocal}

If any non-periodic minimal perfectly balanced set of order $2pq$ can always be combined with another one to form a $2p$-periodic self-complementary set, the reciprocal is also true. Therefore, any $2p$-periodic self-complementary set of order $2pq$ can always be decomposed into two distinct non-periodic minimal perfectly balanced sets.

\begin{theorem} \label{thm:Decomposition_C}
Let $C\in\Z[C_{2pq}]$ be a $2p$-periodic self-complementary set. Then, there exist $A,B \in \Gamma_{2pq}$ such that:
$$C = A + B.$$
\end{theorem}

\begin{proof}
It suffices to show that $C$ can be written up to rotation in the form $C = C_q (P_1 + \delta P_2)$, where $(P_1,P_2)$ is a nontrivial partition of $C_p$. With this expression, by setting $(Q_1,Q_2)$ as a nontrivial partition of $C_q$, we then have:
$$ \begin{aligned}
C &= C_q (P_1 + \delta P_2) \\
 &= (Q_1 + Q_2) (P_1 + \delta P_2) \\
 &= (P_1 Q_1 + \delta P_2 Q_2) + (P_1 Q_2 + \delta P_2 Q_1) \\
 &= A + B
\end{aligned} $$
where $A = P_1 Q_1 + \delta P_2 Q_2$ and $B = P_1 Q_2 + \delta P_2 Q_1$ are two non-periodic minimal perfectly balanced sets by Theorem~\ref{thm:cas2pq} which proves the Theorem.

We start by defining $P_1$ and $P_2$. Since $\Z_{2pq}\simeq \Z_{q} \times \Z_{2p}$, we have $\Z_{2pq}= C_q (C_p + \delta C_p)$. We then define:$ P_1 = \{ x \in C_p | x \in C \} $ and $ P_2 = \{ x \in C_p | \delta x \in C \}.$
By construction, we have $C = C_q (P_1 + \delta P_2)$. It remains to show that $P_1$ and $P_2$ form a nontrivial partition of~$C_p$.

We first show that $P_1$ and $P_2$ are nonempty. By contradiction, if $P_2$ is empty, then $|C| = q \times |P_1|$. However, $C$ is self-complementary, hence $|C| = pq$, which implies that $|P_1| = p$ and therefore $P_1 = C_p$ which leads to $C = C_pC_q$. This implies that $C$ is $2$-periodic which is a contradiction because $C$ is $2p$-periodic. Therefore, $P_2$ is nonempty and the same argument shows that $P_1$ is also nonempty.

We then show that $P_1 + P_2 = C_p$.  Since $C$ is self-complementary, we have:
$$ \begin{aligned}
C_{2pq} &= C + C^c \\
&= C + \delta C \\
&= C_q(P_1+\delta P_2) + C_q(\delta P_1+ P_2) \\
&= C_q(P_1+P_2 + \delta (P_1+ P_2)) \\
&=  C_q C_2 (P_1+ P_2)
\end{aligned} $$
However, $\Z_{2pq}\simeq \Z_{q} \times \Z_{2} \times \Z_{p}$, so $C_{2pq} = C_q C_2 C_p$, which shows that $P_1 + P_2 = C_p$.

Finally, we show that $P_1$ and $P_2$ have no common elements. Because $C$ is self-complementary, we have $|C| = pq$. However, $|C| = q \times (|P_1| + |P_2|)$, hence $|P_1| + |P_2| = p$. Since $P_1 + P_2 = C_p$, we have $|P_1| + |P_2| \geq p$, with equality if and only if they have no elements in common. Therefore, we can conclude that $P_1$ and $P_2$ form a nontrivial partition of $C_p$, which completes the proof. 
\end{proof}

While Figure~\ref{fig:Association_30_1} displays three decompositions into dual perfectly balanced sets of the same $6$-periodic set, Figure~\ref{fig:Association_30_2} illustrates only one decomposition for each $10$-periodic set. The following proposition determines the number of decompositions of a set $C$ that is $2p$-periodic and self-complementary into two dual elements of $\Gamma_{2pq}$.

\begin{proposition} \label{prop:Number_decompositions}
Let $C\in\Z[C_{2pq}]$ be a $2p$-periodic self-complementary set. The number of decompositions of $C$ into two elements of $\Gamma_{2pq}$ up to rotation is: $$\frac{2^{q-1}-1}{q}.$$
\end{proposition}

\begin{proof}
By the proof of Theorem~\ref{thm:Decomposition_C}, there exists a nontrivial partition $(P_1,P_2)$ of $C_p$ such that $C=C_q(P_1+\delta P_2)$. Therefore, every decomposition of $C$ into two dual elements $A,B \in \Gamma_{2pq}$ is uniquely determined by a nontrivial partition $(Q_1,Q_2)$ of $C_q$. Since $C_q$ has $2^q$ partitions into two subsets, and two of them are trivial, there are $2^q-2$ nontrivial partitions into two subsets. Because the sets $A$ and $B$ are non-periodic and $C$ is $2p$-periodic, it remains to divide by $q$, which is the number of rotations required to obtain $C$ using the same $A$ and $B$ up to rotation. Moreover, since exchanging the roles of $A$ and $B$ does not produce a new decomposition, we then divide by 2 which leads to the formula of the proposition.
\end{proof}

\begin{remark}
It is also possible to determine the number of subsets $C$ of $\Z_{2pq}$ that are $2p$-periodic and self-complementary. Since $C$ is determined by a nontrivial partition $(P_1,P_2)$ of $C_p$, that corresponds to counting the nontrivial partitions of $C_p$ modulo rotations and the symmetry obtained by exchanging $P_1$ and $P_2$. Therefore, by an argument similar to that of the previous proof, there are $(2^{p-1}-1)/p$ distinct $2p$-periodic self-complementary subsets of $\Z_{2pq}$ up to rotation. Since there are $(2^{q-1}-1)/q$ distinct decompositions of each such set $C$ into two dual elements $A,B \in \Gamma_{2pq}$ (Proposition~\ref{prop:Number_decompositions}), multiplying these two quantities gives the number of dual associations forming $2p$-periodic self-complementary sets. Finally, since each dual association corresponds to two non-periodic minimal perfectly balanced sets $A$ and $B$, multiplying by $2$ yields the number of elements of $\Gamma_{2pq}$, in agreement with the formula established in Proposition~\ref{prop:Number_of_PBS}.
\end{remark}

While Proposition~\ref{proposition:distribution_6p} explains the distribution and cardinality of non-periodic perfectly balanced sets of order $6p$, the following proposition describes the dual associations of these sets of order $6p$. More precisely, it characterizes the two possible dual associations of the decomposition of $C$: those that are $6$-periodic and those that are $2p$-periodic. First, if $C$ is $6$-periodic and self-complementary, then it is unique up to rotation and is obtained by concatenating the pattern $111000$. Second, if $C$ is $2p$-periodic and self-complementary, then $C$ admits a unique decomposition into two dual non-periodic minimal perfectly balanced sets up to rotation.

\begin{proposition} \label{prop:C_6p_decomposition}
Let $C\in\Z[C_{6p}]$ be a self-complementary set. Up to rotation, we have:
\begin{itemize}
    \item If $C$ is $6$-periodic, then $C$ consists of $p$ concatenations of the pattern 111000.
    \item If $C$ is $2p$-periodic, then $C$ admits a unique decomposition into two elements of $\Gamma_{6p}$.
\end{itemize}
\end{proposition}

\begin{proof}
First, suppose that $C$ is $6$-periodic. Then $C=C_p(Q_1+\delta Q_2)$ with $(Q_1,Q_2)$ a nontrivial partition of $C_3$. Therefore, up to rotation and permutation of $Q_1$ and $Q_2$, we can assume that $Q_1=1+\eta$ and $Q_2=\eta^2$. Expressing them in terms of $z$, we have $\eta=z^{2p}$ and $\delta\eta^2=z^{3p}z^{4p}=z^{7p}$, which leads to:
$$ C=C_p(1+z^{2p}+z^{7p}). $$
Since $C_pz^6=C_p$, the exponents inside the parentheses can be reduced modulo $6$. Moreover, because $p$ is a prime greater than 3, we have $p\equiv\pm1\pmod{6}$. Therefore:
$$
C=
\begin{cases}
C_p(1+z+z^2), & \text{if } p\equiv1\pmod{6},\\
C_p(1+z^{-1}+z^{-2}), & \text{if } p\equiv-1\pmod{6}.
\end{cases}
$$
In both cases, $C$ is equal, up to rotation, to $p$ concatenations of the pattern 111000.

Second, suppose that $C$ is $2p$-periodic. Applying Proposition~\ref{prop:Number_decompositions} with $q=3$, there is therefore $(2^{3-1}-1)/3=1$ dual association of two non-periodic minimal perfectly balanced sets that form $C$ up to rotation.
\end{proof}

\begin{figure}[h]
     \centering
     \vspace{0.3cm}
     \includegraphics[width=0.31\textwidth]{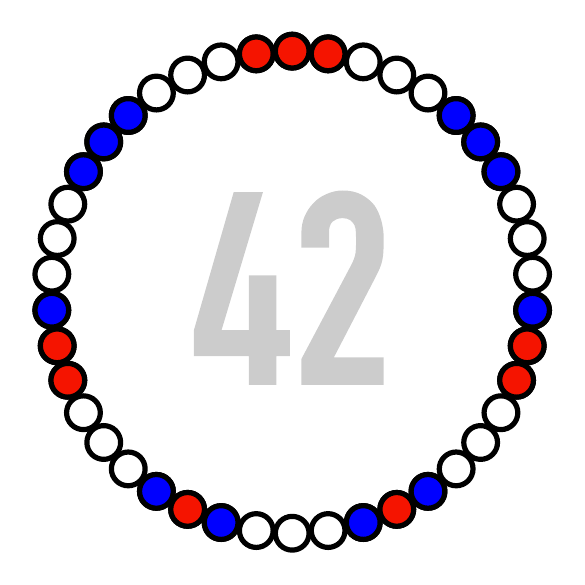} \hfill
     \includegraphics[width=0.31\textwidth]{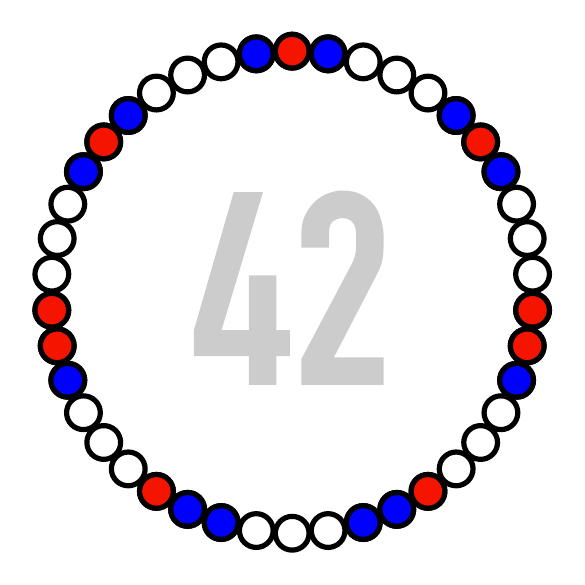} \hfill
     \includegraphics[width=0.31\textwidth]{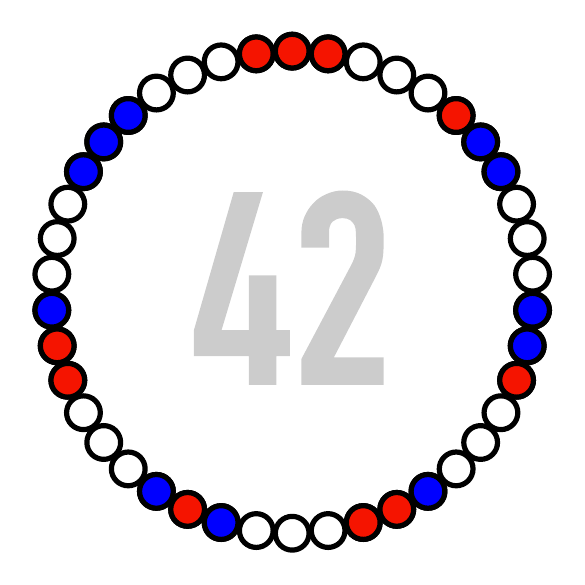}
     \caption{Dual association of non-periodic minimal perfectly balanced sets of $\Z_{42}$ resulting in the same 6-periodic self-complementary set which is a concatenation of the 
     pattern 111000.}
     \label{fig:N_42_Associations}
\end{figure}

Proposition~\ref{prop:C_6p_decomposition} explains, in particular, the six different decompositions in the case of order $6p$ with $p=5$, which are illustrated in Figure~\ref{fig:Association_30_1} and Figure~\ref{fig:Association_30_2}. Firstly, if $C$ is 6-periodic, it always results in 5 concatenations of the pattern 111000 as illustrated in Figure~\ref{fig:Association_30_1}. Secondly, if $C$ is $10$-periodic, then it admits a unique decomposition into two dual perfectly balanced sets, as illustrated by the three different sets $C$ of Figure~\ref{fig:Association_30_2}. An example of the case of order $6p$ with $p=7$ is given in Figure~\ref{fig:N_42_Associations} where the resulting set $C$ is 6-periodic and therefore it always corresponds to 7 concatenations of the pattern 111000. Note that the last resulting set (on the right of Figure~\ref{fig:N_42_Associations}) is composed of perfectly balanced sets which are not symmetric.

\section{Discrete Fourier Transform of Perfectly Balanced Sets} \label{sec:Fourier}

In this section, we present additional theoretical results concerning perfectly balanced sets and their connection with the \textit{discrete Fourier transform} (DFT). This follows recent research on applications of the DFT to music theory, where this approach has been used to study and analyze musical structures defined as subsets of~$\Z_n$~(\citet{amiot2016music}). In particular, we prove some results that we conjectured in a recent paper on the behavior of the DFT modulus of dual perfectly balanced sets~(\citet{lascabettes2026perfectly}).


\subsection{Interpretation of the Discrete Fourier Transform Coefficients of Perfectly Balanced Sets} \label{subsec:Fourier_Interpretation}

Let us first recall the definition of the DFT applied to a subset of $\Z_n$.

\begin{definition}
Let $A \subset \Z_n$. The coefficients $\widehat{A}(k)$ of the \textit{discrete Fourier transform} of $A$ are defined by:
$$ \forall k \in \Z_n, \qquad \widehat{A}(k)=\sum_{a\in A} e^{-2\pi i k a /n}. $$
\end{definition}

The first two coefficients of the DFT of a subset $A$ of $\Z_n$ can be easily interpreted through the following lemma.

\begin{lemma} \label{lem:DFT_1}
Let $A \subset \Z_n$. Then:
$$\widehat{A}(0)= |A| \qquad \text{and} \qquad  \widehat{A}(1)=0 \iff A \text{ is } \textnormal{perfectly balanced}.$$
\end{lemma}

\begin{figure}[h]
    \centering
    \includegraphics[width=0.9\textwidth]{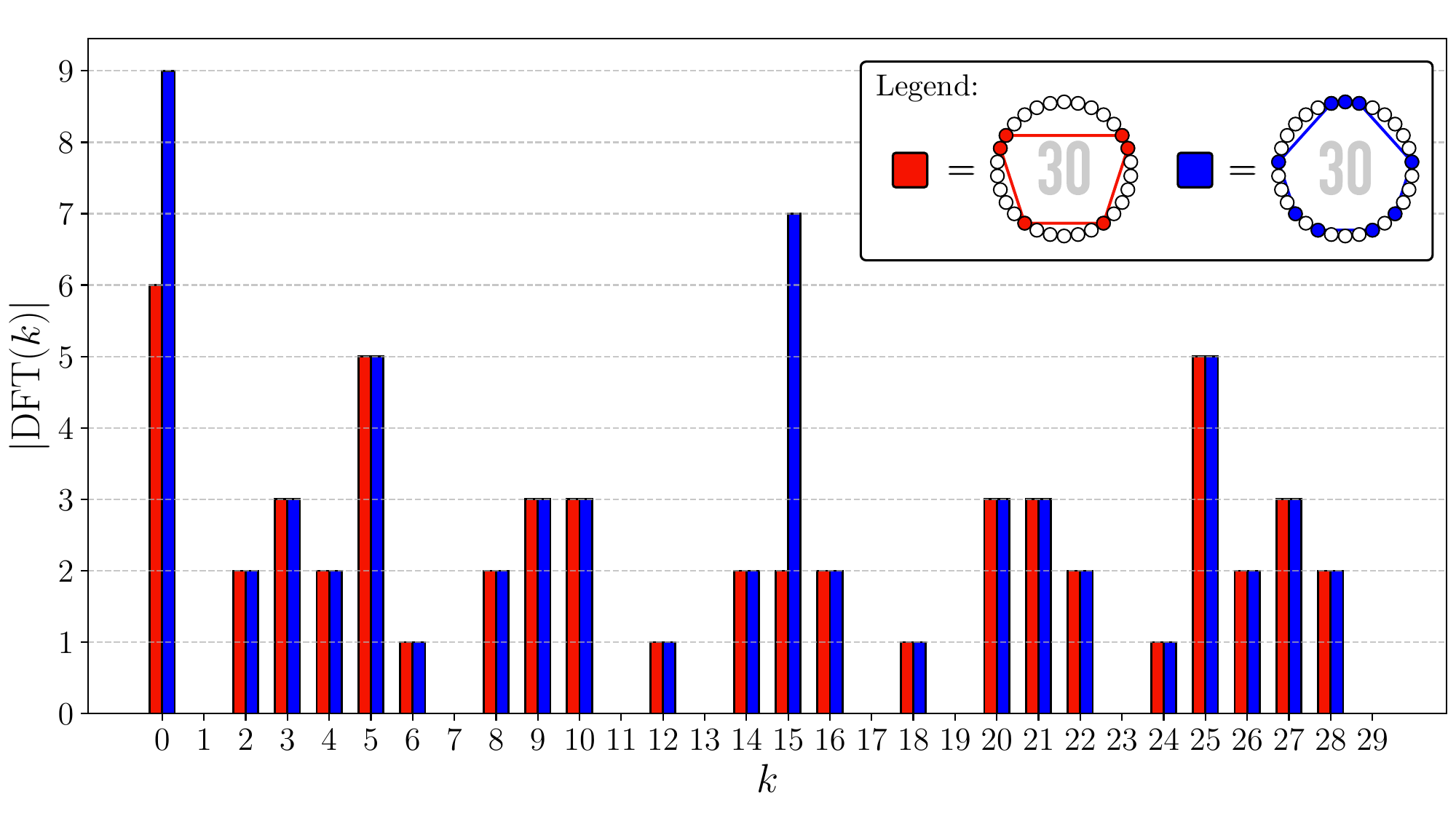}
    \captionsetup{margin=25pt}
    \caption{The coefficients of the modulus of the DFT of two dual perfectly balanced sets of $\Z_{30}$ (the coefficients are identical except at $0$ and $n/2$).}
    \label{fig:DFT}
\end{figure}

Figure~\ref{fig:DFT} displays the modulus of the DFT coefficients of two non-periodic minimal perfectly balanced sets of~$\Z_{30}$. The first set (shown in blue) has a cardinality of 6, and the second (shown in orange) has a cardinality of 9, as can be seen in the figure from the value of the $0^{th}$ coefficient. Moreover, the figure indicates that these two sets are perfectly balanced, since the value of the $1^{st}$ coefficient is indeed zero. It is then possible to interpret other coefficients of the DFT. For example, the coefficients that are coprime to $n$ are also equal to zero for perfectly balanced sets (Theorem 3.11, p.~59, \citet{amiot2016music}). 

\begin{lemma} \label{lem:DFT_2}
Let $A \subset \Z_n$. Then:
$$\forall p \in \Z_n^\times, \quad \widehat{A}(p)=0 \iff A \text{ is } \textnormal{perfectly balanced}.$$
\end{lemma}

In the case of $\Z_{30}$, these are the coefficients with indices 1, 7, 11, 13, 17, 19, 23, and 29, which can be seen in Figure~\ref{fig:DFT} as having zero values due to the perfect balance property. Furthermore, when $n$ is even, $\widehat{A}(n/2)$ is equal to the difference between the numbers of even and odd elements of $A$.

\begin{lemma} \label{lem:DFT_3}
Let $A \subset \Z_n$ with $n$ even. Then:
$$\widehat{A}(n/2)= \sum_{a\in A} (-1)^a.$$
\end{lemma}


\subsection{Relationship Between Discrete Fourier Transform Coefficients of Dual Perfectly Balanced Sets} \label{subsec:Fourier_Relationship}

We now demonstrate a remarkable property concerning the coefficients of the DFT modulus of dual perfectly balanced sets of $\Z_n$ with $n=2pq$: the coefficients are identical everywhere except at $0$ and $n/2$. Figure~\ref{fig:DFT} illustrates this phenomenon for two dual perfectly balanced sets of $\Z_{30}$ where the coefficients coincide except for the $0^{th}$ and $15^{th}$ coefficients. If they are not equal for the $0^{th}$ and $n/2^{th}$ coefficients, this is because these two coefficients are related to the cardinality of the sets (Lemma~\ref{lem:DFT_1} and Lemma~\ref{lem:DFT_3}), and dual perfectly balanced sets do not have the same cardinality. The following result shows that all remaining coefficients always coincide and proves a conjecture about the semi-homometric property of dual perfectly balanced sets (Section 4, p.~259, \citet{lascabettes2026perfectly}).

\begin{proposition}
Let $A,B\subset \Z_n$ be two dual perfectly balanced sets with $n=2pq$. Then:
$$\forall k \in  \Z_n \setminus\{0, n/2\}, \quad \lvert \widehat{A}(k) \rvert = \lvert \widehat{B}(k) \rvert.$$
\end{proposition}

\begin{proof}
By Lemma~\ref{lem:DFT_2}, if $k$ is coprime to $n$, then $\widehat{A}(k)=0$ and $\widehat{B}(k)=0$ since $A$ and $B$ are perfectly balanced, therefore $|\widehat{A}(k)|=|\widehat{B}(k)|$. It remains to consider the case where $k$ is not coprime to $n=2pq$, that is, when $2$, $p$, or $q$ divides $k$. According to Theorem~\ref{thm:cas2pq} and Theorem~\ref{thm:Pairwise_Association}, there exist nontrivial partitions $(P_1, P_2)$ of $C_p$ and $(Q_1, Q_2)$ of $C_q$ such that $A$ and $B$ can be written as:
$$ A = P_1 Q_1 + \delta P_2 Q_2 \quad \text{and} \quad B = P_1 Q_2 + \delta P_2 Q_1.$$
Hence $A+B = C_q(P_1+\delta P_2)$. If $q$ does not divide $k$, then $\widehat{C_q}(k)=0$, and therefore $\widehat{A}(k)+\widehat{B}(k)=0$. It follows that $|\widehat{A}(k)|=|\widehat{B}(k)|$. By a similar argument, we also have $A+\delta B = C_p(Q_1+\delta Q_2)$. In this case, if $p$ does not divide $k$, then $\widehat{A}(k)+\widehat{\delta B}(k)=0$ and because $|\widehat{\delta B}(k)|=|\widehat{B}(k)|$ we have $|\widehat{A}(k)|=|\widehat{B}(k)|$ as well. This completes the proof, since the only coefficients for which both $p$ and $q$ divide $k$ are the $0^{th}$ and $pq^{th}$ coefficients.
\end{proof}


\section{Conclusion} \label{sec:Conclusion}

In this article, we focus on three main aspects of perfectly balanced musical sets of order $2pq$. First, we describe the algebraic construction of such sets by providing a formula that generates all non-periodic minimal perfectly balanced sets of order $2pq$. Moreover, we provide different decompositions into regular polygons and the value of the highest and lowest cardinalities for a given order. Second, we study the combinatorial enumeration of such sets by giving their exact number and analyzing their distribution according to the cardinality. Third, we introduce an original property concerning the duality of perfectly balanced sets of order $2pq$: any such set can be combined with another one to produce a self-complementary $2p$-periodic set and the reciprocal also holds. In addition, we use the discrete Fourier transform to characterize perfectly balanced sets and show that the DFT coefficients of dual perfectly balanced sets are always related. 

\begin{figure}[h]
    \centering
    \includegraphics[width=\textwidth]{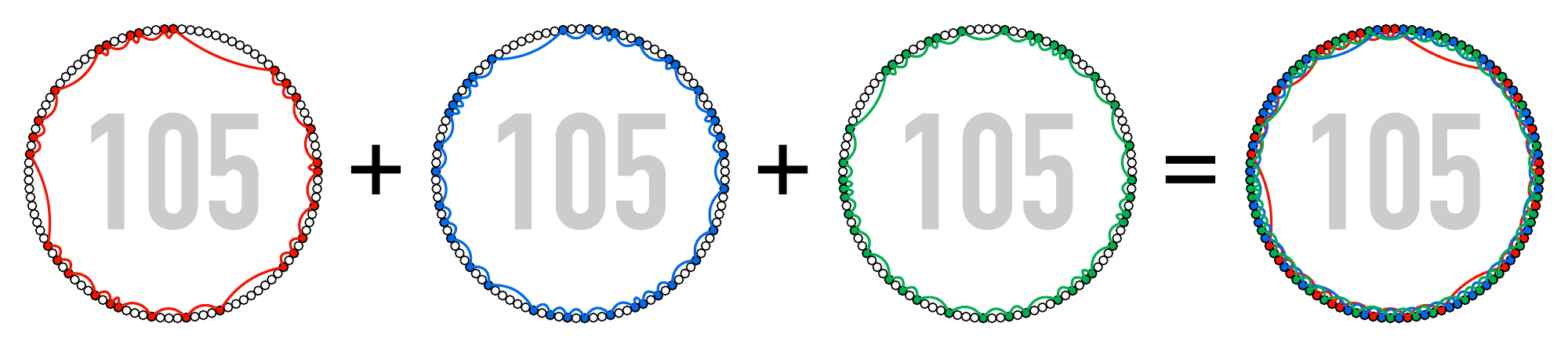}
    \caption{An example of three non-periodic minimal perfectly balanced sets whose union is equal to $\Z_{105}$. Their cardinalities are respectively equal to 29, 37 and 39. This triple association is fundamentally different from the dual association since $n\neq 2pq$ ($105=3 \times 5 \times 7$).}
    \label{fig:N105}
\end{figure}

Although the duality is clear in the case of order $2pq$, its generalization for perfectly balanced sets of a given order remains an open problem. In particular, the association can be of a different nature since the sets are not necessarily pairwise related. For example, Figure~\ref{fig:N105} shows a triple association of perfectly balanced sets of order $3pq$ leading to a complete tiling of the corresponding cyclic group $\Z_{105}$ (since $p=5$ and $q=7$). Note that, in this case, the cardinalities of the three perfectly balanced sets are different (they are respectively equal to 29, 37 and 39). This clearly shows that the association principle, in the general case, is fundamentally different and can be therefore an interesting theoretical problem for future research. Moreover, this also suggests that perfect balance is deeply related to tiling constructions which constitutes an active research field in mathematics and music (\citet{andreatta2009tiling}). Future works may also consider the extension of the perfect balance property to the case of rational coefficients which have a natural application to music since they can be interpreted as note intensity or duration. Finally, all these theoretical results might be implemented into computer-aided environments devoted to the circular representation of musical structures (\citet{milne2016xronomorph,lascabettes2025rhythm}) and therefore be used by music theorists and composers.

As we suggested in the very beginning of this paper, music-theoretical problems can be source of inspiration for advanced mathematical investigations leading to interesting results in both mathematics and music. The perfect balance property is therefore a good example of such a \textit{mathemusical} dynamics, as it can be formalized by using an algebraic framework which leads to original musical objects and new mathematical results. In conclusion, this interdisciplinary approach opens fascinating perspectives towards a mathematical research inspired by music-theoretical constructions.

\section*{CRediT authorship contribution statement}

\textbf{Paul Lascabettes:} Conceptualization, Methodology, Formal analysis, Investigation, Writing -- original draft, Writing -- review \& editing, Visualization.
\textbf{Sidney Ali:} Methodology, Formal analysis, Investigation, Writing -- original draft, Writing -- review \& editing.
\textbf{Alexandre Popoff:} Formal analysis, Investigation, Writing -- review \& editing.
\textbf{Moreno Andreatta:} Methodology, Formal analysis, Writing -- original draft, Writing -- review \& editing, Supervision.

\section*{Declaration of generative AI use}

The authors declare that no AI tools were used in the development of the mathematical results presented in this paper.


\bibliographystyle{plainnat} 
\bibliography{references}

\end{document}